\documentclass[12pt]{article}
\usepackage{amsfonts, amssymb, amsmath}
\usepackage{enumerate}
\usepackage{esint}

\newtheorem{thm}{Theorem}[section]

\newtheorem{lem}[thm]{Lemma}
\newtheorem{prop}[thm]{Proposition}
\newtheorem{defn}[thm]{Definition}
\newtheorem{rem}[thm]{Remark}
\newcommand{\HH}{H_{at}^{1,\vc}(\mu)}
\newcommand{\Hp}{H_{at}^{1,p}(\mu)}
\newcommand{\MS}{M^\sharp}

\newcommand{\f}{\frac}

\newcommand{\al}{\alpha}

\newcommand{\be}{\beta}
\newcommand{\la}{\lambda}
\newcommand{\ep}{\epsilon}
\newcommand{\de}{\delta}
\newcommand{\wt}{\widetilde}

\newcommand{\vc}{\infty}
\newcommand{\dx}{d\mu(x)}
\newcommand{\dy}{d\mu(y)}
\newcommand{\dz}{d\mu(z)}
\newcommand{\dw}{d\mu(w)}

\begin{document}

\title{ Hardy spaces, Regularized BMO spaces and the boundedness of Calder\'on-Zygmund operators on non-homogeneous spaces}
\author{The Anh Bui\thanks{The Anh Bui was supported by a Macquarie University scholarship} \and
Xuan Thinh Duong
\thanks{Xuan Thinh Duong was supported by a research grant from Macquarie University  \newline
{\it {\rm 2010} Mathematics Subject Classification:} 42B20, 42B35.
\newline
{\it Key words:} non-homogeneous spaces, Hardy spaces, BMO,
Calder\'on - Zygmund operator.}}

\date{}


\maketitle
\begin{abstract}
One defines a non-homogeneous space $(X, \mu)$ as a metric space equipped
with a non-doubling measure $\mu$   so that the volume of the ball
with center $x$, radius $r$ has an upper bound of the form $r^n$ for some $n> 0$.
The aim of this paper is to study the boundedness of
Calder\'on-Zygmund singular integral operators $T$
on various function spaces on $(X, \mu)$ such as the Hardy spaces, the
$L^p$ spaces and the regularized BMO spaces. This article
thus extends  the work of X. Tolsa \cite{T1} on the
non-homogeneous space $(\mathbb R^n, \mu)$   to the setting of a general
non-homogeneous space $(X, \mu)$. Our framework of the non-homogeneous
space $(X, \mu)$ is similar to that of \cite{H} and
we are able to obtain quite a few properties similar to those of Calder\'on-Zygmund
operators on doubling spaces such as the weak type $(1,1)$ estimate,
boundedness from Hardy space into $L^1$, boundedness
from $L^{\infty}$ into the regularized BMO and an interpolation theorem. Furthermore, we prove
that the dual space of the Hardy space is the regularized BMO space,
obtain a Calder\'on-Zygmund decomposition
on the non-homogeneous space $(X, \mu)$ and use this decomposition to show the boundedness of
the maximal operators in the form of Cotlar inequality as well as the boundedness of commutators
of Calder\'on-Zygmund operators and BMO functions.
 \end{abstract}

 \tableofcontents
\section{Introduction}

In the last few decades, Calder\'on-Zygmund theory of singular integrals has played a central part of
modern harmonic analysis with lots of extensive applications to other fields of mathematics. This theory has established
criteria for singular integral operators to be bounded on various function spaces including
$L^p$ spaces, $1 < p < \infty$, Hardy spaces, BMO spaces and Besov spaces.
\\

One of the main features of  the standard Calder\'on-Zygmund singular integral theory
is the requirement that the underlying spaces or domains to possess the doubling (volume) property.
Recall that
a space $X$ equipped with a distance $d$ and a measure $\mu$ is said to have the doubling
property if there exists a constant $C$ such that for all $x \in X$ and all $r > 0$,
$$\mu (B(x,2r))  \le C \mu (B(x,r)) $$
where $B(x,r)$ denotes the ball with center $x$ and radius $r > 0$.\\


In the last ten years or so, there has been
substantial progress in obtaining boundedness
of  singular integrals acting  on spaces without the doubling property. Many features of
the standard Calder\'on-Zygmund singular integral theory was extended to spaces with a mild volume growth condition
in place of doubling property through the works of Nazarov, Treil, Volberg, Tolsa
and others. See, for example \cite{NTV1}, \cite{NTV2}, \cite{NTV3}, \cite{T1} and \cite{T2}. These breakthroughs disproved the
long held belief of the decades of 70's and 80's that the doubling property is indispensable in
the theory of Calder\'on-Zygmund singular integrals and lead to more powerful techniques and
estimates in harmonic analysis.
\\

Let  $X$ be a metric space equipped with a measure
$\mu$, possibly non-doubling,  satisfying
$$
\mu(B(x,r))\leq Cr^n
$$
for some positive constants $C$, $n$ and all $r > 0$. We will call such a space $(X, \mu)$
a non-homogeneous space.
Let $T$ be a   Calder\'on-Zygmund operator acting on a non-homogeneous space $X$,
i.e. the associated kernel
of $T$ satisfies appropriate bounds and has H\"older continuity
(for the precise definition, see Section 2.1). Assume that $T$ is bounded on $L^2(X)$,
 then it is shown in  \cite{NTV2} that  the Calder\'on-Zygmund operator  $T$ is of weak type $(1,1)$,
hence by interpolation, is bounded on $L^p(X), 1<p<\vc$.
See also \cite{T2}.\\

Hardy spaces and BMO spaces on a non-homogeneous space $X$ were studied by a number of authors,
for example \cite{T1}, \cite{MMNO}, \cite{H}. In \cite{MMNO}, the authors studied the spaces
BMO$(\mu)$ and $H^1_{at}(\mu)$ on $\mathbb{R}^n$ (with BMO$(\mu)$ space being defined via the standard bounded
oscillations and the Hardy space $H^1_{at}(\mu)$ being defined by an atomic decomposition) for a non doubling measure $\mu$ and
 showed some standard properties of these spaces such as
 the John-Nirenberg inequality, an interpolation
theorem between BMO$(\mu)$ and $H^1_{at}(\mu)$, and
BMO$(\mu)$ being the dual space of $H^1_{at}(\mu)$.
However,  it is shown by Verdera \cite{V}  that an $L^2$ bounded Calder\'on-Zygmund operator may
be unbounded from $L^\vc(\mu)$ into BMO$(\mu)$ as well as
from $H^1_{at}(\mu)$ into $L^1(\mu)$. This shows the need to introduce variants of the BMO spaces characterized by bounded oscillation
estimates so that the Calder\'on-Zygmund operators are bounded from $L^\vc(\mu)$ into these variants of BMO spaces.\\

In \cite{T1}, the author introduced the RBMO space, a variant of the space BMO,
on the non-homogeneous space $( \mathbb{R}^n, \mu)$ which retains some of the
properties of the standard BMO such as the John-Nirenberg
inequality. See Section 3 for the definition of RBMO spaces.
While  Calder\'on-Zygmund operators might not be bounded from
$L^\vc(\mathbb{R}^n, \mu)$ into BMO$(\mathbb{R}^n, \mu)$, they are bounded from
$L^\vc(\mathbb{R}^n, \mu)$ into RBMO$(\mathbb{R}^n, \mu)$, \cite{T1}.\\

Recently,   Hyt\"onen studied the RBMO spaces on non-homogeneous spaces $(X, \mu)$
(instead of $(\mathbb{R}^n, \mu)$) \cite{H}. He proved that the space
RBMO$(\mu)$ on $X$ still satisfies John-Nirenberg inequality. However,
 the boundedness of Calder\'on-Zygmund operators from $L^\vc(\mu)$ into RBMO$(\mu)$
and a number of other properties are still open questions for the setting of general
non-homogeneous spaces $(X, \mu)$.\\

  In this article, our aim  is to conduct an extensive study on
the RBMO spaces  on general non-homogeneous spaces. More specifically, for a
non-homogeneous space $(X,\mu)$ equipped with a measure $\mu$
which is dominated by some doubling measure (the same setting as in \cite{H}),
we are able to prove the following new results:
\begin{enumerate}[(i)]
\item An $L^2$ bounded Calder\'on-Zygmund
operator is bounded  from $L^\vc (\mu)$ into
the RBMO space, see Theorem \ref{boundedofCZOonRBMO}.
\item The dual space of the atomic Hardy spaces is shown to be the RBMO space.
We also show that an $L^2$ bounded Calder\'on-Zygmund
operator is bounded from the atomic Hardy space  $H^1_{at}(\mu)$ into
$L^1(\mu)$, see Theorems \ref{dualHardy} and \ref{boundedofCZOonHardyspace}.
\item An interpolation theorem between the RBMO space and the Hardy space  $H^1_{at}(\mu)$: if an
operator is bounded from $H^1_{at}(\mu)$ into
$L^1(\mu)$ and from $L^\vc (\mu)$ into
the RBMO space, then the operator is bounded on $L^p(\mu)$ for all $1 < p < \infty$,
see Theorem \ref{interpolatioHardyRBMO}.
\item A Cotlar type  inequality for Calder\'on-Zygmund operators which gives
the boundedness of several maximal operators associated with $T$, see Theorem \ref{Cotlar1}.
\item The boundedness of commutators of Calder\'on-Zygmund
operators and RBMO functions on $L^p$ spaces, see Theorem \ref{boundednessofcommutators}.
\item A Calder\'on-Zygmund decomposition using a variant of Vitali covering lemma,
see Theorem \ref{CZdecomposition}, and the weak type $(1,1)$ of an $L^2$ bounded
Calder\'on-Zygmund operator, see Theorem \ref{weak1-1}.
\end{enumerate}

We remark that, while this manuscript was in finishing touch, we
learned that similar results concerning the Hardy spaces as in (ii) have
been obtained independently in \cite{HYY}. \\

We now give a brief comment about some techniques used  in this paper.
 In addition to using some ideas and techniques
in \cite{T1}, we obtain certain key estimates through careful
investigation of the family of doubling balls in a non-homogeneous space $(X, \mu)$.
 Let us recall that the
main techniques used in \cite{T1} rely on the Besicovitch covering lemma and the
construction of the $(\al,\be)$-doubling balls in $\mathbb{R}^n$. However, the
Besicovitch covering lemma is only applicable to $\mathbb R^n$ and it is not
applicable in the setting of  general non-homogeneous spaces. In the general setting,
one can construction the $(\al,\be)$-doubling
balls by using a covering lemma in \cite{He} in place of the Besicovitch covering lemma.
In \cite{H},  the author  used this substitution to
obtain the John-Nirenberg inequality for the RBMO spaces.
However, it seems that to obtain further results similar to the standard theory as in the case of doubling spaces,
 more refined techniques are needed.\\

 An important technical detail in this paper  is our
construction of  the three consecutive $(\al, \be)$-doubling balls
(see, Proposition \ref{rem1}) which we employ successfully to obtain
the important characterizations (\ref{cond1-RBMOdefn3}) and
(\ref{cond2-RBMOdefn3}), similar to those  in \cite[Lemma 2.10]{T1}. By using these
 three consecutive $(\al, \be)$-doubling balls,  we show the
boundedness of Calder\'on-Zygmund operators from $L^\vc (X, \mu)$ into
the space RBMO  (see Theorem \ref{boundedofCZOonRBMO}) as well as
an interpolation theorem of RBMO spaces (see Theorem \ref{interpolationthm}).\\


\bigskip

{\it Acknowledgement: } The second named author would like to thank El Maati Ouhabaz
for helpful discussion.\\


\section{Non-homogeneous spaces, families of doubling balls and singular integrals}

\subsection{Non-homogeneous spaces and families of doubling balls}

In this paper, for the sake of simplicity we always assume that
$(X,d)$ is a metric space. With
minor modifications,  similar  results  hold when $X$ is a
quasi-metric space.\\

{\bf Geometrically doubling regular metric spaces.} We adopt the
definition that the space $(X, d)$ is geometrically doubling if
there exists a number $N\in \mathbb{N}$ such that every open ball
$B(x, r) = \{y \in X : d(y, x) < r\}$ can be covered by at most $N$
balls of radius $r/2$. Our using of  this somewhat non-standard name
is to differentiate this property from other types of doubling
properties. If there is no specification, the ball $B$ means the
ball center $x_B$ with radius $r_B$. Also, we set $n = \log_2 N$,
which can be viewed as (an upper bound for) a geometric dimension of
the space. Let us recall the following well-known lemma. See, for
example \cite{H}.
\begin{lem}\label{lem1.1}
In a geometrically doubling regular metric space, a ball $B(x, r)$ can
contain the centers $x_i$ of at most $N\alpha^{-n}$ disjoint balls
$B(x_i, \alpha r)$ for  any $\alpha\in (0, 1]$.
\end{lem}

{\bf Upper doubling measures.} A measure $\mu$ in the metric space $(X,\mu)$ is
said to be an upper doubling measure if there exists a dominating
function $\lambda$ with the following properties:
\begin{enumerate}[(i)]
\item $\lambda:X\times (0,\vc)\mapsto (0,\vc)$;
\item for any fixed $x\in X$, $r\mapsto \lambda(x,r)$ is increasing;
\item there exists a constant $C_\lambda > 0$ such that
$
\lambda(x,2r)\leq C_\lambda \lambda(x,r)
$
for all $x\in X$, $r >0$;
\item  the  inequality
$
\mu(x,r) := \mu(B(x,r)) \leq \lambda(x,r)
$
holds for all $x\in X$, $r >0$;
\item and $\lambda(x,r)\approx \lambda(y,r)$ for all $r>0 , \ x, y\in X$
and $d(x,y)\leq r$.
\end{enumerate}
We note that in \cite{H}, the condition (v) is not assumed.
\begin{lem}\label{coveringlemma}
Every family of balls $\{B_i\}_{i\in F}$ of uniformly bounded
diameter in a metric space $X$ contains a disjoint sub-family
$\{B_i\}_{i\in E}$ with $E\subset F$ such that
$$
\cup_{i\in F}B_i\subset \cup_{i\in E}5B_i.
$$
\end{lem}
For a proof of Lemma \ref{coveringlemma}, see \cite{He}.\\

{\bf Assumptions: Throughout the paper, we always assume that $(X, \mu)$ is a
geometrically doubling regular metric space and the measure $\mu$
is an upper doubling measure.}\\

We adopt the following definition as in \cite{T1}. For $\al,\be>1$, a ball $B\subset X$ is called $(\al,\be)$-doubling
if $\mu(\alpha)\leq \be\mu(B)$. The following result states the
existence of plenty of  doubling balls with small radii and with large radii.\\

\begin{lem}[\cite{H}]\label{lem1}
The following statements hold:
\begin{enumerate}[(i)]
\item If $\be>C_\lambda^{\log_2\al}$, then for any ball $B\subset X$ there exists $j\in \mathbb{N}$ such that $\al^jB$ is $(\al,\be)$-doubling.
\item If $\be>\al^n$  where $n$ is the doubling order of $\lambda$, then for any ball $B\subset X$ there exists $j\in \mathbb{N}$ such that $\al^{-j}B$ is $(\al,\be)$-doubling.
\end{enumerate}
\end{lem}

\bigskip

Our following result which shows the existence of three
consecutive $(\al,\be)$ doubling balls will play an important role
in this paper.\\
\begin{prop}\label{rem1}
If $B$ is a $(\alpha^3,\beta)$ doubling ball ($\alpha>1$), then $B,
\alpha B$ and $\alpha^2 B$ are three consecutive $(\alpha,\beta)$
doubling balls.
\end{prop}
\emph{Proof:} The proof of Proposition \ref{rem1} is simple,
hence we omit the details here.

\medskip

For any two balls $B\subset Q$, we defined
\begin{equation}\label{coefficientK}
K_{B, Q}=1+\int_{r_{B}\leq d(x,x_{B})\leq
r_{Q}}\f{1}{\lambda(x_{B},d(x,x_{B}))}\dx.
\end{equation}
This definition is a variant of the definition  in \cite[pp.94-95]{T1}.
Similarly to the results \cite[Lemma 2.1]{T1} we have the following
properties:
 \begin{lem}\label{lem2}
\begin{enumerate}[(i)]
\item If $Q\subset R\subset S$ are balls in $X$, then $$\max\{K_{Q,R},K_{R,S}\}\leq K_{Q,S}\leq C(K_{Q,R}+K_{R,S}).$$
\item If $Q\subset R$ are compatible size, then $K_{Q,R}\leq C.$
\item If $\al Q,\ldots \al^{N-1}Q$ are non $(\al,\be)$-doubling balls $(\be>C_\lambda^{\log_2\al})$ then $K_{Q,\al^N Q}\leq C.$
\end{enumerate}
\end{lem}
The proof of Lemma \ref{lem2} is not difficult, hence we omit the
details here.

\medskip

As in \cite{T1}, for two balls $B\subset Q$ we can define the
coefficient $K'_{B, Q}$ as follows: let $N_{B, Q}$ be the smallest
integer satisfying $6^{N_{B, Q}}r_B\geq r_{Q}$, then we set
$$
K'_{B,Q}:=1+\sum_{k=1}^{N_{B,Q}}\f{\mu(6^kB)}{\lambda(x_B, 6^kr_B)}.
$$
In the case that $\lambda(x,ar)=a^m\lambda(x,r)$ for all $x\in X$ and $a, r>0$, it is
not difficult to show that $K_{B,Q}\approx K'_{B,Q}$. However, in
general, we only have $K_{B,Q}\leq C K'_{B,Q}$.
\subsection{Calder\'on-Zygmund operators}

A kernel $K(\cdot,\cdot)\in L^1_{{\rm loc}}(X\times X\backslash
\{(x,y): x=y\})$ is called a Calder\'on-Zygmund kernel if
\begin{enumerate}[(i)]
\item \begin{equation}\label{cond1ofC-Z}|K(x,y)|\leq C\min\Big\{\f{1}{\lambda(x,d(x,y))},\f{1}{\lambda(y,d(x,y))}\Big\}.
\end{equation}
\item There exists $0<\delta\leq 1$ such that
\begin{equation}\label{cond2ofC-Z}
|K(x,y)-K(x',y)|+|K(y,x)-K(y,x')|\leq
C\f{d(x,x')^\delta}{d(x,y)^\delta \lambda(x,d(x,y))}
\end{equation}
if $d(x,x')\leq Cd(x,y).$
\end{enumerate}

A linear operator $T$ is called a Calder\'on-Zygmund operator with
kernel $K(\cdot, \cdot)$ satisfying (\ref{cond1ofC-Z}) and
(\ref{cond2ofC-Z}) if for all $f\in L^\vc (\mu)$ with bounded
support and $x\notin {supp}f$,
$$
Tf(x)=\int_X K(x,y)f(y)\dy.
$$
The maximal operator $T_*$ associated with the Calder\'on-Zygmund
operator $T$ is defined by
$$
T_*f(x)=\sup_{\ep>0}|T_\ep f(x)|,
$$
where $T_\ep f(x)=\int_{d(x,y)\geq \epsilon} K(x,y)f(y)\dy$.\\

We would like to give an example for the operator whose the
associated kernel satisfies the conditions (\ref{cond1ofC-Z}) and
(\ref{cond2ofC-Z}). As in \cite{H}, we consider  Bergman-type operators which are
studied by Volberg and Wick. In \cite{VW}, the authors obtained a
characterization of measures $\mu$ in the unit ball
$\mathbb{B}_{2n}$ of $\mathbb{C}^n$ for which the analytic
Besov-Sobolev space $B^\sigma_2 (\mathbb{B}^{2n})$ embeds
continuously into $L^2(\mu)$. Their proof goes through a new $T1$
theorem for what they call �Bergman-type� operators. Let us describe
the situation of  this application. The measures $\mu$ in \cite{VW}
satisfy the upper power bound $\mu(B(x, r)) \leq r^m$, except
possibly when $B(x, r) \subset H$, where $H$ is a fixed open set.
However, in the exceptional case there holds $r\leq \delta(x) :=
d(x,H^c)$, and hence $$ \mu(B(x, r)) \leq \lim_{\epsilon\rightarrow
0} B(x, \delta(x) + \epsilon)\leq \lim_{\epsilon\rightarrow 0}
(\delta(x) + \epsilon)^m=\delta^m.$$ Thus we find that their
measures are actually upper doubling with
$$\mu (B(x, r))\leq \max\{\delta(x)^m, r^m) =: \lambda(x, r).$$ It is
not difficult to show that $\lambda(\cdot,\cdot)$ satisfies the
conditions (i)-(v) in definition of upper doubling measures.\\

In \cite{VW}, as a main application concerning the
Besov-Sobolev spaces, the authors introduced the operator associated
to the kernel
\begin{equation}\label{kerneldef}
K(x, y) = (1 - \overline{x}\cdot y)^{-m},
\end{equation}
for $x, y \in
\overline{\mathbb{B}}_{2n}\subset \mathbb{C}^n$.
 Here $\overline{x}$
stands for the componentwise complex conjugation, and $\overline{x}\cdot y$ designates
the usual dot product of $n$-vectors $\overline{x}$ and $y$. Moreover, one equips
$\overline{\mathbb{B}}_{2n}$ with the regular quasi-distance, see
\cite[Lemma 2.6]{Tch},
$$
d(x, y) :=\Big||x| - |y|\Big| +\Big| 1 - \f{\overline{x}\cdot y}{
|x| |y|}\Big|.
$$
Finally, the set $H$ related to the exceptional balls is now the
open unit ball $\overline{\mathbb{B}}_{2n}$. It was proved in \cite{HM}
that the kernel $K(x,y)$ defined by (\ref{kerneldef}) satisfies
(\ref{cond1ofC-Z}) and (\ref{cond2ofC-Z}).

\section{The RBMO spaces}
\subsection{Definition of RBMO$(\mu)$}
The RBMO (Regularized BMO) space was introduced by Tolsa for $(\mathbb R^n, \mu)$ in \cite{T1}
and it was adopted by T. Hyt\"onen for general non-homogeneous space $(X, \mu)$ in \cite{H}.
\begin{defn}
Fix a parameter $\rho>1$. A function $f\in L^1_{{\rm loc}}(\mu)$ is
said to be in the space {\rm RBMO}$(\mu)$ if there exists a number
$C$, and for every ball $B$, a number $f_B$ such that
\begin{equation}\label{cond1-RBMOdefn1}
\f{1}{\mu(\rho B)}\int_B|f(x)-f_B|\dx\leq A
\end{equation}
and, for any two balls $B$ and $B_1$ such that  $B\subset B_1$,
\begin{equation}\label{cond2-RBMOdefn1}
|f_B-f_{B_1}|\leq C K_{B,B_1}.
\end{equation}
The infimum of the values $C$  in (\ref{cond2-RBMOdefn1}) is taken to be the {\rm RBMO} norm of $f$ and denoted
by $\|f\|_{{\rm RBMO}(\mu)}$.
\end{defn}
The {\rm RBMO} norm  $\|\cdot\|_{{\rm RBMO}(\mu)}$ is independent of
$\rho>1$. Moreover the John-Nirenberg inequality holds for {\rm
RBMO}$(X)$. More precisely, we have the following result (see Corollary
6.3 in \cite{H}).
\begin{prop}
For any $\rho >1$ and $p\in [1,\vc)$, there exists a constant $C$ so
that for every $f\in  {\rm RBMO}(\mu)$ and every ball $B_0$,
$$
\Big(\f{1}{\mu(\rho
B_0)}\int_{B_0}|f(x)-f_{B_0}|^p\dx\Big)^{1/p}\leq C\|f\|_{{\rm
RBMO}(\mu)}.
$$
\end{prop}

\subsection{Some characterizations of RBMO$(\mu)$}
In the rest of paper, unless $\al$ and $\be$
are specified otherwise, by an $(\alpha,\beta)$ doubling ball we mean a $(6, \beta_0)$-doubling with a fixed number $\be_0>\max\{C_\lambda^{3\log_2 6},6^{3n}\}$.\\

Given a ball $B\subset X$, let $N$ be the smallest non-negative integer such that $\widetilde{B}=6^NB$ is doubling. Such a ball $\widetilde{B}$ exists due to Lemma \ref{lem2}.\\

Let $\rho>1$ be some fixed constant. We say that $f\in L^1_{{\rm
loc}}(\mu)$ is in {\rm RBMO}$(\mu)$ if there exists some constant
$C>0$ such that for any ball $Q$
\begin{equation}\label{cond1-RBMOdefn2}
\f{1}{\mu(\rho B)}\int_B |f(x)-m_{\wt{B}}f|\dx\leq C
\end{equation}
and
\begin{equation}\label{cond2-RBMOdefn2}
|m_Qf-m_Rf|\leq CK_{Q,R}, \ \ \text{for any two doubling balls
$Q\subset R$},
\end{equation}
here $m_Bf$ is the mean value of $f$ over the ball $B$. Then we take
$$\|f\|_{*}:=\inf\{C: \text{(\ref{cond1-RBMOdefn2}) and (\ref{cond2-RBMOdefn2}) hold}\} .$$\\
By the same proof as in Lemma 2.8 of \cite{T1}, we have the
following result.
\begin{prop}
For a fixed $\rho>1$, the norms $\|\cdot\|_*$ and $\|\cdot\|_{{\rm
RBMO}(\mu)}$ are equivalent.
\end{prop}

We now extend certain characterizations of ${\rm RBMO}(\mu)$ in \cite{T1}
in the case of $(\mathbb{R}^n, \mu)$ to the case of non-homogeneous spaces $(X, \mu)$.
 In the case of $\mathbb{R}^n$, Besicovitch covering lemma was used but this lemma is not applicable in our setting.
 We overcome this problem by using the three consecutive doubling balls in Proposition \ref{rem1}.

\begin{prop}\label{otherdefnofRBMO}
For $f\in L^1_{{\rm loc}}(\mu)$, the following are equivalent:
\begin{enumerate}[(a)]
\item $f\in {\rm RBMO}(\mu)$.
\item There exists some constant $C_b$ such that for any ball $B$
\begin{equation}\label{cond1-RBMOdefn3}
\f{1}{\mu(6 B)}\int_B |f(x)-m_Bf|\dx\leq C_b
\end{equation}
and
\begin{equation}\label{cond2-RBMOdefn3}
|m_Qf-m_Rf|\leq C_bK_{Q,R}\Big(\f{\mu(6 Q)}{\mu(Q)}+\f{\mu(6
R)}{\mu(R)}\Big), \ \ \text{for any two balls $Q\subset R$}.
\end{equation}
\item There exists some constant $C_c$ such that for any doubling ball $B$
\begin{equation}\label{cond1-RBMOdefn4}
\f{1}{\mu(B)}\int_B |f(x)-m_Bf|\dx\leq C_c
\end{equation}
and
\begin{equation}\label{cond2-RBMOdefn4}
|m_Qf-m_Rf|\leq C_cK_{Q,R}, \ \ \text{for any two doubling balls
$Q\subset R$}.
\end{equation}
\end{enumerate}
Moreover, the best constants $C_b$ and $C_c$ are comparable to the
{\rm RBMO}$(\mu)$ norm of $f$.
\end{prop}

\emph{Proof:} $(a)\rightarrow (b):$ If $f\in {\rm
RBMO}(\mu)$, then (\ref{cond1-RBMOdefn3}) and
(\ref{cond2-RBMOdefn3}) hold for $C_b=C\|f\|_*$ for some constant $C$.
Indeed, for any ball $B$ we have
$$
|m_Bf -m_{\wt{B}}f|\leq m_Q(|f-m_{\wt{B}}f|)\leq
\|f\|_*\f{\mu(6Q)}{Q}.
$$
Therefore,
\begin{equation}\label{eq1}
\f{1}{\mu(6 B)}\int_B |f(x)-m_Bf|\dx\leq \f{1}{\mu(6 B)}\int_B
(|f-m_{\wt{B}}f|+|m_Bf -m_{\wt{B}}f|)\leq 2\|f\|_*.
\end{equation}
On the other hand, for any two balls $Q\subset R$, one has
$$
|m_Qf-m_Rf|\leq |m_Qf -m_{\wt{Q}}f|+|m_{\wt{Q}}f
-m_{\wt{R}}f|+|m_{R}f-m_{\wt{R}}f|.
$$
Applying (\ref{eq1}) for the first and the third terms, we have
$$
|m_Qf -m_{\wt{Q}}f|+|m_{R}f-m_{\wt{R}}f|\leq
\|f\|_*\Big(\f{\mu(6Q)}{Q}+\f{\mu(6R)}{R}\Big).
$$
We can follow the argument in \cite{T1} to obtain the estimate for
the second term. Let us remark that for any two balls $Q\subset R$
such that $\wt{Q}\subset \wt{R}$, it follows from (\ref{cond2-RBMOdefn2})  that
$$
|m_{\wt{Q}}f -m_{\wt{R}}f|\leq \|f\|_* K_{\wt{Q},\wt{R}}.
$$
By Lemma \ref{lem2}, we have
$$
K_{\wt{Q},\wt{R}}\leq C(K_{Q,\wt{Q}} +K_{Q, R}+K_{R,\wt{R}})\leq
C(C_1 +K_{Q, R}+C_2)\leq CK_{Q, R}.
$$
In general, $Q\subset R$ does not imply $\wt{Q}\subset \wt{R}$. We consider two cases:\\
{\it Case 1:} If $r_{\wt{Q}}\geq r_{\wt{R}}$, then $\wt{Q}\subset
3\wt{R}$. Setting $R_0= \wt{3\wt{R}}$, then it follows from Lemma
\ref{lem2} and (\ref{cond2-RBMOdefn2}) that
\begin{equation*}
\begin{aligned}
|m_{\wt{Q}}f -m_{\wt{R}}f|&\leq |m_{\wt{Q}}f -m_{R_0}f|+m_{R_0}f -m_{\wt{R}}f\\
&\leq (K_{\wt{Q},R_0}+K_{\wt{R},R_0})\|f\|_*.
\end{aligned}
\end{equation*}
For the term $K_{\wt{Q},R_0}$ we have
\begin{equation*}
\begin{aligned}
K_{\wt{Q},R_0}&\leq CK_{Q,R_0}\\
&\leq C(K_{Q,R}+K_{R,R_0})\\
&\leq C(K_{Q,R}+K_{R,\wt{R}} +K_{\wt{R}, 3\wt{R}}+K_{3\wt{R},R_0})\\
&\leq CK_{Q,R}.
\end{aligned}
\end{equation*}
The remaining term $K_{\wt{R},R_0}$ is dominated by
$$
C(K_{\wt{R},3\wt{R}}+K_{3\wt{R},R_0})\leq CK_{Q,R}.
$$
So in this case, we obtain $|m_{\wt{Q}}f -m_{\wt{R}}f|\leq CK_{Q,R}\|f\|_*$.\\

{\it Case 2:} If $r_{\wt{R}}<r_{\wt{Q}}$, then $\wt{R}\subset 6^2
\wt{Q}$. Obviously, we can find some $m\geq 1$ such that
$r_{\wt{R}}\geq \f{r_{5^mQ}}{25}$ and $\wt{R}\subset 6^mQ\subset
6^2\widetilde{Q}$. Therefore, $\wt{R}$ and $5^mQ$ are comparable
sizes. This implies $K_{\wt{R},5^mQ}\leq C$. Setting $Q_0=\wt{6^2
\wt{Q}}$ we have
\begin{equation*}
\begin{aligned}
|m_{\wt{Q}}f -m_{\wt{R}}f|&\leq |m_{\wt{Q}}f -m_{Q_0}f| + |m_{Q_0}f -m_{\wt{R}}f|\\
&\leq (K_{\wt{Q},Q_0}+K_{\wt{R},Q_0})\|f\|_*.
\end{aligned}
\end{equation*}
Let us estimate $K_{\wt{Q},Q_0}$. We have
$$
K_{\wt{Q},Q_0}\leq C(K_{\wt{Q},6^2 \wt{Q}}+K_{6^2\wt{Q},Q_0})\leq CK_{Q,R}.
$$
For the term $K_{\wt{R},Q_0}$, one has
\begin{equation*}
\begin{aligned}
K_{\wt{R},Q_0}&\leq
C(K_{\wt{R},5^mQ}+K_{5^mQ,6^2\wt{Q}}+K_{6^2\wt{Q},Q_0})\\
&\leq
C(K_{\wt{R},5^mQ}+K_{Q,6^2\wt{Q}}+K_{6^2\wt{Q},Q_0})\\
&\leq CK_{Q,R}.
\end{aligned}
\end{equation*}
Therefore, in this case we also obtain $|m_{\wt{Q}}f -m_{\wt{R}}f|\leq CK_{Q,R}\|f\|_*$.\\

$(b) \rightarrow (c)$: the proof of this implication it easy and hence we omit the detail here.\\

$(c) \rightarrow (a)$: Let $B$ be some ball. We need
to show that (\ref{cond1-RBMOdefn2}) holds for $\rho=6$. For
any $x\in B$, there exists some $(6^3,\beta_0)$-doubling ball
centered $x$ with radius $r_{6^{-2j}B}$ for some $j\in \mathbb{N}$.
We denote by $B_x$  the biggest ball satisfying these properties. Let
us recall that by Proposition \ref{rem1}, the balls $B_x, 6B_x$ and $6^2B_x$
are three $(6,\beta_0)$-doubling balls. Moreover, by Lemma \ref{lem2} we
have
$$
|m_{6B_x}f-m_{\wt{B}}f|\leq
|m_{6B_x}f-m_{B_x}f|+|m_{B_x}f-m_{\wt{B}}f|\leq CC_c.
$$
By Lemma \ref{coveringlemma}, we can pick a disjoint subcollection
$B_{x_i}$, $i\in I$, such that $B\subset \cup_{i\in
I}5B_{x_i}\subset \cup_{i\in I}6B_{x_i}$. Thus, we have
\begin{equation*}
\begin{aligned}
\int_B|f-m_{\wt{B}}f|d\mu&\leq \sum_{i\in I} \int_{B_{x_i}}|f-m_{\wt{B}}f|d\mu\\
&\leq \sum_{i\in I} \int_{B_{x_i}}(|f-m_{6B_{x_i}}f|+|m_{6B_{x_i}}f-m_{\wt{B}}f|)d\mu\\
&\leq \sum_{i\in I} CC_c\mu(6B_{x_i})\\
&\leq \sum_{i\in I} C\be_0C_c\mu(B_{x_i})\\
&\leq C\be_0C_c\mu(6B).\\
\end{aligned}
\end{equation*}
This completes our proof.

\section{Interpolation results}

\subsection{The sharp maximal operator}
Adapting an idea in \cite{T1}, we define the sharp maximal
operator as follows:
\begin{equation}\label{sharpoperator}
M^\sharp f(x)=\sup_{B\ni x}\f{1}{\mu(6
B)}\int_B|f-m_{\wt{B}}f|d\mu+\sup\limits_{(Q,R)\in
\Delta_x}\f{|m_Qf-m_Rf|}{K_{Q,R}},
\end{equation}
here $\Delta_x:=\{(Q,R): x\in Q\subset R \ \text{and} \  Q,R: {\rm doubling}\}$.\\

Note that in our sharp maximal operator, the term $\mu(6B)$ was chosen with the fixed constant $6$ throughout the paper.
It is clear that
$$
f\in {\rm RBMO}(\mu) \Leftrightarrow \MS f\in L^\vc(\mu).
$$
We define, for $\rho\geq 1$, the non-centered maximal operator
$M_{(\rho)}$ by setting
$$
M_{(\rho)}f(x)=\sup_{x\in Q}\f{1}{\mu(\rho Q)}\int_Q|f|d\mu.
$$
It was proved that $M_{(\rho)}$ is of weak type $(1,1)$ for
$\rho\geq 5$ and hence $M_{(\rho)}$ is bounded on $L^p(\mu)$ for all
$p\in (1,\vc]$, see \cite[Proposition 3.5]{H}. When $\rho=1$, we write $Mf$ instead
of $M_{(1)}f$. From the boundedness of $M_{(\rho)}$ for $\rho\geq
5$, the non-centered doubling maximal operator  is defined by
$$
Nf(x)=\sup_{x\in Q: \ {\rm doubling}}\f{1}{\mu(Q)}\int_Q|f|d\mu
$$
where the supremum is taken over all $(6, \beta_0)$ doubling balls, is of weak type $(1,1)$ and hence bounded on $L^p(\mu)$ for all $p\in
(1,\vc]$.\\
Note that it is not difficult to show that
$$
M^\sharp f(x)\leq M_{(6)}f(x) + 3Nf(x)
$$
for all $x\in X$. Therefore the operator $M^\sharp$ is of type weak $(1,1)$ and bounded on $L^p(\mu)$ for all $1<p<\vc$.
\begin{lem}\label{sharpmaximalestimate1}
For $f\in L^1_{{\rm loc}}(\mu)$, we have $$\MS
|f|(x) \leq 5\be_0 \MS f(x).$$
\end{lem}
The proof is similar to that of Remark 6.1 in \cite{T1}.\\

We now show that the non-centered doubling maximal operator is
dominated by the sharp maximal operator in the following theorem.
Although, some estimates are inspired from \cite[Theorem 6.2]{T1},
there are some main differences in our proof. More specifically, the
three consecutive doubling balls
argument will be used to replace
the Besicovitch covering lemma.

\begin{thm}\label{sharpmaximalestimate2}
Let $f\in L^1_{\rm loc}(\mu)$ with the extra condition $\int f d\mu =0$ if
$\|\mu\|:=\mu(X)<\vc$. Assume that for some $p$, $1<p<\vc$, $\inf\{1, Nf\}\in L^p(\mu)$.
Then we have
$$
\|Nf\|_{L^p(\mu)}\leq C\|\MS f\|_{L^p(\mu)}.
$$
\end{thm}
\emph{Proof:} We assume that $\|\mu\|=\vc$. The proof for $\|\mu\|<\vc$
is similar. By
standard argument, it suffices to prove the following $\lambda$-good
inequality: for some fixed $\nu<1$ and all $\epsilon>0$ there exists
some $\de >0$ such that for any $\lambda>0$ we have
\begin{equation}\label{lambda goodinequatlity}
\mu\{x: Nf(x)>(1+\epsilon)\lambda, \MS f(x)\leq \delta\lambda\}\leq
\nu\mu\{x: Nf(x)>\lambda\}.
\end{equation}

Setting $E_\lambda=\{x: Nf(x)>(1+\epsilon)\lambda, \MS f(x)\leq
\delta\lambda\}$ and $\Omega_\lambda=\{x: Nf(x)>\lambda\}$, for
$f\in L^p(\mu)$. For each $x\in E_\lambda$, we can choose the
doubling ball $Q_x$ containing $x$ satisfying that
$m_{Q_x}|f|>(1+\ep/2)\la$ and if $Q$ is any doubling ball containing
$x$ with $r_Q> 2r_{Q_x}$ then $m_Q|f|\leq (1+\ep/2)\la$. Such a ball
$Q_x$ exists due to $f\in L^p(\mu)$.

\medskip

Let $R_x$ be the ball centered $x$ with radius $6r_{Q_x}$ and $S_x$
be the smallest $(6^3,\be_0)$-doubling ball in the form $6^{3j}R_x$.
Then, by Proposition \ref{rem1}, $S_x, 6S_x$ and $6^2S_x$ are three
$(6, \be_0)$-doubling balls. Moreover, one has
$$
K_{Q_x, 6S_x}\leq C(K_{Q_x, R_x}+K_{R_x, S_x}+K_{S_x, 6S_x})\leq C.
$$
Therefore, it follows from Lemma \ref{sharpmaximalestimate1} that
$$
|m_{Q_x}|f|-m_{6S_x}|f||\leq K_{Q_x, 6S_x}\MS |f|(x)\leq C\beta_0
\MS f(x)\leq C\be_0\de\la.
$$
This implies that for sufficiently small $\de$ we have
$$
m_{6S_x}|f|>\la
$$
and hence $6S_x\subset \Omega_\la$.

\medskip

Note that by Lemma \ref{coveringlemma}, we can pick a
disjoint collection $\{S_{x_i}\}_{i\in I}$ with $x_i\in E_\la$ and
$E_\lambda\subset \cup_{i\in I}5S_{x_i}\subset \cup_{i\in
I}6S_{x_i}$. Setting $W_{x_i}=6S_{x_i}$, we will show
that
\begin{equation}\label{equ1}
\mu(6S_{x_i}\cap E_\la)\leq C\f{\nu}{\be_0}\mu(6S_{x_i})
\end{equation}
for all $i\in I$.

\medskip

Once (\ref{equ1}) is proved, (\ref{lambda goodinequatlity}) follows readily. Indeed, from (\ref{equ1}) we have
$$
\mu(E_\la)\leq \sum_{i\in I}\mu(6S_{x_i}\cap E_\la)\leq \sum_{i\in
I}\f{\nu}{\be_0}\mu(6S_{x_i})\leq C\sum_{i\in I}\nu\mu(S_{x_i})\leq
C\nu\mu(\Omega_\la).
$$

Now we show the proof of (\ref{equ1}). Let $y\in
W_{x_i}\cap E_\la$. For any doubling ball $Q\ni y$
satisfying $m_Q|f|>(1+\ep)\la$, it follows that $r_Q\leq
r_{W_{x_i}}/8$. Indeed, if $r_Q>
r_{W_{x_i}}/8$ then we have $Q_{x_i}\subset
W_{x_i}\subset \wt{16Q}$ and
$$
|m_Q|f|-m_{\wt{16Q}}|f||\leq K_{Q,\wt{16Q}}\MS |f|(y)\leq
C\de\la\leq \f{\ep}{2}
$$
for sufficiently small $\de$.
This implies $m_{\wt{16Q}}|f|>(1+\ep/2)\la$ which is a contradiction
to the choice of $Q_{x_i}$. So, $r_Q\leq r_{W_{x_i}}/8$.
This, together with $m_Q|f|>(1+\ep)\la$, imply
$$
N(f\chi_{\f{5}{4}W_{x_i}})(y)>(1+\ep)\la
$$
and
$$
m_{\wt{\f{5}{4}W_{x_i}}}|f|\leq (1+\ep/2)\la \
\text{(since
$r_{\widetilde{\f{5}{4}W_{x_i}}}>2r_{Q_{x_i}}$}).
$$
This yields,
$$
N(\chi_{\f{5}{4}W_{x_i}}|f|-m_{\wt{\f{5}{4}W_{x_i}}}|f|)(y)>\f{\ep}{2}\la.
$$
 Therefore, by using the weak $(1,1)$ boundedness of $N$, we have
 \begin{equation*}
 \begin{aligned}
 \mu(W_{x_i}\cap E_\la)&\leq \mu\{y:N(\chi_{\f{5}{4}W_{x_i}}|f|-m_{\wt{\f{5}{4}W_{x_i}}}|f|)(y)>\f{\ep}{2}\la\}\\
 &\leq \f{C}{\ep\la}\int_{\f{5}{4}W_{x_i}}(|f|-m_{\wt{\f{5}{4}W_{x_i}}}|f|)d\mu\\
 &\leq \f{C}{\ep\la}\mu(\f{15}{2}W_{x_i})\MS |f|(x_i)\\
 &\leq \f{C\de}{\ep}\be_0\mu(6^3S_{x_i})\\
&\leq \f{C\de}{\ep}\be_0\mu(S_{x_i}).
\end{aligned}
\end{equation*}
Thus, (\ref{equ1}) holds provided $\delta <\ep/C\nu\be_0$.\\
For the case $f\notin L^p(\mu)$, we define the sequence of functions
$\{f_k\}, \ k = 1, 2, \cdots$ by setting
\begin{equation*}
f_k(x)=\begin{cases}
f(x), &|f(x)|\leq k,\\
k\f{f(x)}{|f(x)|}, \ & |f(x)|>k.
\end{cases}
\end{equation*}
Then we have $\MS f_k(x)\leq C\MS f(x)$. On the other hand,
$|f_k(x)|\leq k\inf\{1,|f|(x)\}\leq k\inf(1,Nf)(x)$ and so $f_k\in
L^p(\mu)$. Hence,
$$\|Nf_k\|_{L^p(\mu)}\leq C\|\MS f_k\|_{L^p(\mu)}\leq C\|\MS f\|_{L^p(\mu)}.$$
Taking the limit as $k\rightarrow\infty$, we obtain the required result and the proof is completed.

\subsection{An Interpolation Theorem for linear operators}
\begin{thm}\label{interpolationthm}
Let $1<p<\infty$ and let $T$ be a linear operator bounded on $L^p(\mu)$
and from $L^\vc(\mu)$ into ${\rm RBMO(\mu)}$. Then $T$ extends to a
bounded operator on  $L^r(\mu)$ for $p<r<\infty$.
\end{thm}
\emph{Proof:} We consider 2 cases:\\
{\bf Case 1: $\|\mu\|=\vc$:} Since $T$ is bounded on $L^p(\mu)$,
$\MS T$ is sublinear bounded on $L^p(\mu)$ and on $L^\vc(\mu)$.
Therefore, by interpolation, $\MS T$ is bounded on
$L^r(\mu)$ for $p<r<\infty$,
$$
\|\MS Tf\|_{L^r(\mu)}\leq C\|f\|_{L^r(\mu)}.
$$
Assume that $f\in L^r(\mu)$ is supported in compact set. Then $f\in
L^p(\mu)$ and so $Tf\in L^p(\mu)$.  Hence $Nf\in L^p(\mu)$ and
$\inf \{1,Nf \} \in L^r(\mu)$. By invoking Theorem
\ref{sharpmaximalestimate2},
$$
\|Tf\|_{L^r(\mu)}\leq\|\MS Tf\|_{L^r(\mu)}\leq C\|f\|_{L^r(\mu)}.
$$
{\bf Case 2}: Assume that $\|\mu\|<\vc$. For $f\in L^r(\mu)$, set $f=(f-\int f
d\mu)+\int f d\mu=f_1+f_2$. Since $\int f_1 d\mu=0$, we can apply
the same argument as for $\|\mu\|=\vc$. It is not difficult to show
that $\|T1\|_{L^r(\mu)}\leq C\|1\|_{L^r(\mu)}$. This completes the proof.

\section{Atomic Hardy spaces and their dual spaces}

\subsection{The space $H^{1,\vc}_{at}(\mu)$}

For a fixed $\rho>1$, a function $b\in L^1_{loc}(\mu)$ is called
an atomic block if
\begin{enumerate}[(i)]
\item there exists some ball $B$ such that ${\rm supp} b\subset B$;
\item $\int b d \mu=0;$
\item there are functions $a_j$ supported on cubes $B_j\subset B$
and numbers $\lambda_j\in \mathbb{R}$ such that
\begin{equation}\label{eq1-defnHardy}
b=\sum_{j=1}^\vc \lambda_j a_j,
\end{equation}
where the sum converges in $L^1(\mu)$, and
$ \|a_j\|_{L^\vc(\mu)}\leq (\mu(\rho B_j)K_{B_j,B})^{-1}
$
and the constant $K_{B_j,B}$ being  given in the paragraph before Lemma \ref{lem2}.
\end{enumerate}

We denote $|b|_{H^{1,\vc}_{at}(\mu)}=\sum_{j=1}^\vc
|\lambda_j|$. We say that $f\in H^{1,\vc}_{at}(\mu)$ if there are
atomic blocks $b_i$ such that
\begin{equation}\label{eq2-defnHardy}
f=\sum_{i=1}^\vc b_i
\end{equation}
with $\sum_{i=1}^\vc |b_i|_{H^{1,\vc}_{at}(\mu)}<\vc$. The
$H^{1,\vc}_{at}(\mu)$ norm of $f$ is defined by
$$
\|f\|_{H^{1,\vc}_{at}(\mu)}:=\inf\sum_{i=1}^\vc
|b_i|_{H^{1,\vc}_{at}(\mu)}
$$
where the infimum is taken over all the possible decompositions of
$f$ in atomic blocks.
\medskip

We have the following basic properties of $H_{at}^{1,\vc}(\mu)$.
\begin{prop}\label{prop1HH}
\begin{enumerate}[(a)]
\item $\HH$ is a Banach space.
\item $\HH\subset L^1(\mu)$ and $\|f\|_{L^1(\mu)}\leq \|f\|_{\HH}$.
\item The space $\HH$ is independent of the constant $\rho$ when
$\rho>1$.
\end{enumerate}
\end{prop}
\emph{Proof:} The proofs of (a) and (b) are standard and we omit the details here.

\medskip
The proof of (c): Given $\rho_1>\rho_2>0$, it is clear that
$H^{1,\vc}_{at,\rho_1}\subset H^{1,\vc}_{at,\rho_2}$ with
$\|f\|_{H^{1,\vc}_{at,\rho_2}}\leq \|f\|_{H^{1,\vc}_{at,\rho_1}}$.
Conversely, if $b=\sum_{i=1}^\vc\lambda_i a_i$ is an atomic block
with supp $a_i\subset B_i\subset B$ in $H^{1,\vc}_{at,\rho_1}$, then
by Lemma \ref{lem1.1} we can cover each $B_i$ by $N
\Big[\f{\rho_1}{\rho_2}\Big]^n$ balls, says $\{B_{ik}\}$, with the
same radius $\f{\rho_2}{\rho_1}r_B$. Therefore, we can decompose
$a_i:=\sum_{k}a_{ik}$ where
$a_{ik}:=a_i\f{\chi_{B_{ik}}}{\sum_{j}\chi_{B_{ij}}}$. It is not
difficult to verify that $b$ is also an atomic block in
$H^{1,\vc}_{at,\rho_2}$. This completes our proof.\\

We now show that the space ${\rm RBMO}(\mu)$ is embedded in the dual space of $\HH.$

\begin{lem}\label{lem1-hardyspace}  We have
$${\rm RBMO}(\mu)\subset \HH^*.$$
That is, for $g\in {\rm RBMO}(\mu)$, the linear functional
$$
L_g(f)=\int_X  fgd\mu
$$
defines a continuous linear functional $L_g$ over $\HH$ with
$
\|L_g\|_{\HH^*}\leq C\|g\|_{{\rm RBMO}(\mu)}.
$
\end{lem}
\emph{Proof:} Following standard argument, see for example
\cite[p.64]{CW2}, we only need to check that for an atomic
block $b$ and $g\in $RBMO$(\mu)$, we have
$$
\Big|\int bg d\mu\Big|\leq C|b|_{\HH}\|g\|_{{\rm RBMO}(\mu)}.
$$
Assume that supp$b\subset B$ and $b=\sum_{j}^\vc \lambda_j a_j$,
where $a_j$'s are functions satisfying (a) and (b) in the definition of
atomic blocks. If $g\in L^\vc$, by using $\int bd\mu = 0$, we have
\begin{equation}\label{eq1-proofH-BMO}
\Big|\int bg d\mu\Big|=\Big|\int b(g-g_B) d\mu\Big|\leq \sum_{j}^\vc
|\lambda_j| \|a_j\|_{L^\vc(\mu)}\int_{B_i}|g-g_B|d\mu.
\end{equation}
Since $g\in L^\vc(\mu)\subset {\rm RBMO}(\mu)$, we have
\begin{equation*}
\begin{aligned}
\int_{B_i}|g-g_B|d\mu&\leq
\int_{B_i}|g-g_{B_i}|d\mu+\int_{B_i}|g_B-g_{B_i}|d\mu\\
&\leq CK_{B_i,B}\|g\|_{{\rm RBMO}(\mu)}\mu(\rho B_j).
\end{aligned}
\end{equation*}
From (\ref{eq1-proofH-BMO}), we obtain
$$
\Big|\int bg d\mu\Big|\leq C|b|_{\HH}\|g\|_{{\rm RBMO}(\mu)}.
$$
In general case, if $g\in {\rm RBMO}(\mu)$, define
$$
g_N(x):=\begin{cases} f(x), &\ \ |f(x)|<N,\\
 N\f{f(x)}{|f(x)|}, &\ \ |f(x)|\geq N.
\end{cases}
$$
It can be verified that $\|g_N\|_{{\rm RBMO}(\mu)}\leq C \|g\|_{{\rm
RBMO}(\mu)}$. As above, since $g_N\in L^\vc(\mu)$, we have
$$
\Big|\int fg_N d\mu\Big|\leq C\|f\|_{\HH}\|g_N\|_{{\rm
RBMO}(\mu)}\leq C\|f\|_{\HH}\|g\|_{{\rm RBMO}(\mu)}.
$$
Let us denote $L^\vc_{0}:=\{f: f \ \text{in $L^\vc(\mu)$ with
bounded support}\}$ and $D=\HH\cap L^\vc_{0}$. So, the functional
$L_g: f\mapsto \int gf$ is well-defined on $D$ whenever $g\in {\rm
RBMO}(\mu)$ (since $g\in L^1_{loc}(\mu)$). By the dominated
convergence theorem
$$
\lim_{N\rightarrow \vc}\int fg_N d\mu =\int fg d\mu
$$
for all $f\in D$. We claim that $D$ is dense in $\HH$. To verify this claim,
denote by $H^{1,\vc}_{at,fin}(\mu)$ the set of  all elements in $\HH$
where the sums (\ref{eq1-defnHardy}) and (\ref{eq2-defnHardy})
are taken over finite elements. Obviously,
$H^{1,\vc}_{at,fin}(\mu)$ is dense in $\HH$ and each functional
$f\in H^{1,\vc}_{at,fin}(\mu)$ is also in $L^\vc_{0}$. Therefore ,
$L_b$ is a unique extension on $\HH$ and hence
$$
\Big|\int fg d\mu\Big|\leq C\|f\|_{\HH}\|g\|_{{\rm RBMO}(\mu)}.
$$
This completes our proof. \\

The following lemma can be obtained by the same argument as in
\cite[Lemma 4.4]{T1}.
\begin{lem}
If $g\in {\rm RBMO}(\mu)$, we have
$$
\|L_g\|_{\HH}\approx \|g\|_{{\rm RBMO}(\mu)}.
$$
\end{lem}

\subsection{The space $H^{1,p}_{at}(\mu)$}
For a fixed $\rho>1$, a function $b\in L^1_{loc}(\mu)$ is called a
$p$-atomic block, $1<p<\vc$, if
\begin{enumerate}[(i)]
\item there exists some ball $B$ such that ${\rm supp} b\subset B$;
\item $\int b d \mu=0;$
\item there are functions $a_j$ supported on cubes $B_j\subset B$
and numbers $\lambda_j\in \mathbb{R}$ such that
\begin{equation}\label{eq1-defnHardy p}
b=\sum_{j=1}^\vc \lambda_j a_j,
\end{equation}
where the sum converges in $L^1(\mu)$, and
$$
\|a_j\|_{L^p(\mu)}\leq (\mu(\rho B_j))^{1/p-1}K_{B_j,B}^{-1}.
$$
\end{enumerate}
We denote $|b|_{H^{1,p}_{at}(\mu)}=\sum_{j=1}^\vc |\lambda_j|$.
We say that $f\in H^{1,p}_{at}(\mu)$ if there are
$p$-atomic blocks $b_i$ such that
\begin{equation}\label{eq2-defnHardy p}
f=\sum_{i=1}^\vc b_i
\end{equation}
with $\sum_{i=1}^\vc |b_i|_{H^{1,p}_{at}(\mu)}<\vc$. The
$H^{1,p}_{at}(\mu)$ norm of $f$ is defined by
$$
\|f\|_{H^{1,p}_{at}(\mu)}:=\inf\sum_{i=1}^\vc
|b_i|_{H^{1,p}_{at}(\mu)}
$$
where the infimum is taken over all the possible decompositions of
$f$ in $p$-atomic blocks.\\

Similarly to $\HH$, we have the following basic properties of
$H_{at}^{1,p}(\mu)$
\begin{prop}
\begin{enumerate}[(a)]
\item $\Hp$ is a Banach space.
\item $\Hp\subset L^1(\mu)$ and $\|f\|_{L^1(\mu)}\leq \|f\|_{\Hp}$.
\item The space $\Hp$ is independent of the constant $\rho$ when
$\rho>1$.
\end{enumerate}
\end{prop}
The proofs of this proposition is in line with Proposition \ref{prop1HH}, so we omit the details here.

\medskip

\begin{lem}\label{lem1-hardyspace-p-atom} We have
$${\rm RBMO}(\mu)\subset \Hp^*.$$
That is, for $g\in {\rm RBMO}(\mu)$, the linear functional
$$
L_g(f)=\int_X fgd\mu
$$
defines a continuous linear functional $L_g$ over $\Hp$ with
$$
\|L_g\|_{\Hp^*}\leq C\|g\|_{{\rm RBMO}(\mu)}.
$$
\end{lem}
\emph{Proof:} The proof of this lemma is analogous to that of Lemma
\ref{lem1-hardyspace} with minor modifications. We
leave the details to the interested reader.\\

We remark that a main difference between the Hardy space in
Tolsa's setting \cite{T1} and our Hardy space in this article is the sense of convergence
in the atomic decomposition.
This leads to different approaches in proving the inclusions
RBMO$(\mu) \subset \HH^*$ and RBMO$(\mu) \subset \Hp^*$. However,
for the inverse inclusion $\Hp^* \subset$ RBMO$(\mu)$, by a careful
investigation, Tolsa \cite{T1} showed that one only needs to consider the sums in
(\ref{eq1-defnHardy p}) and (\ref{eq2-defnHardy p}) over finite
$p$-atoms and $p$-atomic blocks, hence the sense of
convergence in (\ref{eq1-defnHardy p}) and (\ref{eq2-defnHardy p})
does not matter in both settings. This is the reason why we can use the arguments
in \cite{T1} for our setting with minor modifications to obtain the duality result of
$\HH$ and $\Hp$ as in the next Theorem.

\begin{thm}\label{dualHardy}

For $1<p<\vc, \Hp=\HH$. Also $\Hp^*=\HH^*={\rm RBMO}(\mu)$.
\end{thm}

As explained above, we omit the details of the proof.

\section{Calder\'on-Zygmund decomposition}
\subsection{Calder\'on-Zygmund decomposition}

The following two technical lemmas will be useful for the construction of a
 Calder\'on-Zygmund decomposition on non-homogeneous spaces.

 \begin{lem}\label{lem1-CZ}
 Assume that $Q$, $S$ are  two concentric balls, $Q\subset R$,  such that there are no
$(\alpha, \beta)$-doubling balls with $\beta>C_\lambda^{\log_2\alpha}$ in
the form $\alpha^kQ, k\in \mathbb{N}$ such that
$Q\subset\alpha^kQ\subset R$. Then we have
$$
\int_{R\backslash Q}\f{1}{\lambda(x_Q,d(x_Q,x))}\dx\leq C.
$$
\end{lem}
\emph{Proof:} Let $N$ be the smallest integer such that $R\subset
\alpha^NQ$. Then, $\mu(\alpha^k Q)\geq \beta \mu(\alpha^{k-1} Q)$
for all $k=1, \ldots, N$. Therefore, we have,
\begin{equation*}
\begin{aligned}
\int_{R\backslash Q}&\f{1}{\lambda(x_Q,d(x_Q,x))}\dx\\
&\leq \sum_{k=1}^N \int_{\alpha^{k-1}r_Q\leq d(x,y)\leq \alpha^kr_Q}\f{1}{\lambda(x_Q,d(x_Q,x))}\dx\\
&\leq \sum_{k=1}^N \f{\mu(\alpha^kQ)}{\lambda(x_Q,\alpha^{k-1}r_Q)}\\
&\leq \sum_{k=1}^N
\f{\beta^{N-k}\mu(\alpha^NQ)}{(C_\lambda)^{(N-k)\log_2\alpha}\lambda(x_Q,\alpha^{N}r_Q)}\\
&\leq \sum_{k=1}^N
\Big[\f{\beta}{(C_\lambda)^{\log_2\alpha}}\Big]^{N-k}\\
&\leq \sum_{j=1}^\vc
\Big[\f{\beta}{(C_\lambda)^{\log_2\alpha}}\Big]^{j}\\
&\leq C \ \ \text{(since $\beta>\log_2\alpha$)}.
 \end{aligned}
\end{equation*}
This completes the proof.\\

While the Covering Lemma \ref{coveringlemma} for $(X, \mu)$ can be used to replace the
Besicovich covering lemma for $(\mathbb R^n, \mu)$ in certain estimates,
 the Calder\'on-Zygmund decomposition in  $(X, \mu)$ will need a covering lemma
 which gives the finite overlapping property at all points $x \in X$.  This is given in the next lemma.

 \begin{lem}\label{lem2-CZ}
Every family of balls $\{B_i\}_{i\in F}$ of uniformly bounded
diameter in a metric space $X$ contains a disjoint sub-family
$\{B_i\}_{i\in E}$ with $E\subset F$ such that
\begin{enumerate}[(i)]
\item $ \cup_{i\in F}B_i\subset \cup_{i\in E}6B_i,
$
\item For each $x\in X$, $\sum_{i\in E}\chi_{6B_i}<\vc$.
\end{enumerate}
\end{lem}

We remark that in (ii), the sum $\sum_{i\in E}\chi_{6B_i}<\vc$ at each $x$ but
these sums are not necessarily uniformly bounded on $X$.

\medskip

\emph{Proof:} By Lemma \ref{coveringlemma} we can pick a disjoint
subfamily $\{B_i: B_i=B(x_{B_i},r_{B_i})\}_{i\in E}$ with $E\subset
F$ satisfying (i). Moreover, we can assume that for
$i,j\in E$, neither   $6B_i\subset 6B_j$ nor $6B_j\subset
6B_i$.

\medskip

To prove (ii), we assume in contradiction that there exists some
$x\in X$ such that there exists an infinite family of balls $\{B_i:
i\in I_x \subset E \}$ such that $x\in B_i$ for all $i\in I_x$. We
will show that $\lim\inf_{i\in I_x} r_{B_i}>0$. Otherwise, for any
$\epsilon>0$ there exists $i_{\epsilon}\in I_x$ such that
$r_{B_{i_\epsilon}}<\epsilon$. Therefore, if $B_0$ is any ball in
the family $\{B_i: i\in I_x\}$, there exists $r>0$ such that
$B(x,r)\subset 6B_0$. For $\epsilon=\f{r}{30}$, we have $x\in
6B_{i_\epsilon}$ and $r_{6B_{i_\epsilon}}<\f{r}{4}$. This implies
$6B_{i_\epsilon}\subset 6B_0$ which is a contradiction.

\medskip

Thus $\lim\inf_{i\in I_x} r_{B_i}>0$. This together with the uniform
boundedness of diameter of the family of balls shows that there exist
$m$ and $M>0$ such that $m<r_{B_i}<M$ for all $i\in I_x$. Obviously,
$\cup_{i\in I_x} B(x_{B_i},m)\subset B(x, 2M)$ and the balls
$\{B(x_{B_i},m): i\in I_x\}$ are pairwise disjoint. By Lemma
\ref{lem1.1}, there exists a finite family of balls with radius
$\f{m}{30}$ such that $B(x, 2M)\subset \cup_{i=1}^K B(x_i,
\f{m}{30})$. Therefore, there exist a ball, says $B_k\in \{B(x_i,
\f{m}{30}): i\in 1,\ldots K\}$, and at least two balls $B_1$ and
$B_2$ in $\{B_i: i\in I_x\}$ such that $B_k\cap \f{1}{6}B_1\notin
\emptyset$ and $B_k\cap \f{1}{6}B_1\notin \emptyset$. Since
$\min\{r_{\f{1}{6}B_1}, r_{\f{1}{6}B_1}\}>\f{1}{6}m=5r_{B_k}$, we
have $B_k\subset B_1\cap B_2$. This is a contradiction, because the
family of balls $\{B(x_{B_i},m): i\in I_x\}$ is pairwise disjoint.
Our proof is completed.

\medskip

We now give a Calder\'on-Zygmund
decomposition on a non-homogenous space $(X, \mu)$ which is an extension
of a Calder\'on-Zygmund
decomposition on the non-homogeneous space $(\mathbb R^n, \mu)$ in \cite{T1}.

\begin{thm}\label{CZdecomposition}(Calder\'on-Zygmund decomposition)
Assume $1\leq p<\vc$. For any $f\in L^p(\mu)$ and any $\lambda>0$
(with $\lambda> \beta_0||f||_p/||\mu||$ if $||\mu||<\infty$), the
following statements hold.

\begin{description}
\item (a) There exists a family of finite overlapping balls
$\{6Q_i\}_i$ such that $\{Q_i\}_i$ is a pairwise disjoint family and
\begin{equation}\label{cz1}
\frac{1}{\mu(6^2Q_i)}\int_{Q_i}|f|^p d\mu>\frac{\lambda^p}{\beta_0},
\end{equation}
\begin{equation}\label{cz2}
\frac{1}{\mu(6^2 \eta Q_i)}\int_{\eta Q_i}|f|^pd\mu \leq
\frac{\lambda^p}{\beta_0}, \ \text{for all $\eta>1$},
\end{equation}
\begin{equation}\label{cz3}
|f|\leq \lambda \ \text{a.e. $(\mu)$ on
$X \backslash\bigcup_{i}6Q_i$}.
\end{equation}
\item (b) For each $i$, let $R_i$ be a $(3\times6^2, C_\lambda^{\log_2 3\times 6^2 +1})$-doubling
ball concentric with $Q_i$, with $l(R_i)>6^2l(Q_i)$ and denote
$\omega_i=\frac{\chi_{6Q_i}}{\sum_k\chi_{6Q_k}}$. Then there exists
a family of functions $\varphi_i$ with constant signs and supp
$(\varphi_i)\subset R_i$ satisfying
\begin{equation}\label{cz4}
\int\varphi_id\mu=\int_{6Q_i}f\omega_id\mu,
\end{equation}
\begin{equation}\label{cz5}
\sum_i|\varphi_i|\leq \kappa \lambda,
\end{equation}
(where $\kappa$ is some constant which depends only on $(X, \mu)$), and
\begin{enumerate}[(i)]
\item \begin{equation}\label{cz6} ||\varphi_i||_{\infty}\mu(R_i)\leq
C\int_{X}|w_if|d\mu \ \text{if $p=1$};
\end{equation}
\item  \begin{equation}\label{cz6.1} ||\varphi_i||_{L^p(\mu)}\mu(R_i)^{1/p'}\leq
\f{C}{\lambda^{p-1}}\int_{X}|w_if|^p d\mu \ \text{if $1<p<\vc$}.
\end{equation}
\end{enumerate}
\item (c) For $1<p<\vc$, if $R_i$ is the smallest $(3\times6^2, C_\lambda^{\log_2 3\times 6^2 +1})$-doubling
ball of the family $\{3\times 6^2 Q_i\}_{k\geq 1}$, then
\begin{equation}\label{cz8}
\|b\|_{\Hp}\leq \f{C}{\lambda^{p-1}}\|f\|_{L^p(\mu)}^p
\end{equation}
where $b=\sum_{i}(w_if-\varphi_i))$.
\end{description}
\end{thm}
\emph{Proof:} For the sake of simplicity, we only give the proof for
the case $p=1$ for (a) and (b). When $p>1$, by setting $g=f^p\in
L^1(\mu)$, we can reduce to the problem $p=1$. Then, with a simple
modification,
we will obtain (\ref{cz6.1}) instead of (\ref{cz6}).\\

(a) Set $E:=\{x: |f(x)|>\lambda\}$. For each $x\in E$, there exists
some ball $Q_x$ such that
\begin{equation}\label{cz7}
\frac{1}{\mu(6^2Q_x)}\int_{Q_x}|f|d\mu>\frac{\lambda}{\beta_0}
\end{equation}
and such that if $Q'_{x}$ is centered at $x$ with $l(Q'_x)>l(Q_x)$,
then
\begin{equation*}
\frac{1}{\mu(6^2Q'_x)}\int_{Q'_x}|f|d\mu\leq \frac{\lambda}{\beta_0}
\end{equation*}
Now we can apply Lemma \ref{lem2-CZ} to get a family of balls
$\{Q_i\}_{i}\subset \{Q_x\}_x$ such that
$\sum_{j}\chi_{6Q_j}(x)<\vc$ for all $x\in X$ and (\ref{cz1}), (\ref{cz2}) and (\ref{cz3}) are satisfied.

\medskip

(b) Assume first that the family of balls $\{Q_i\}$ is finite.
Without loss of generality, suppose that $l(R_i)\leq l(R_{i+1})$.
The functions $\varphi$ will be constructed of the form
$\varphi_i=\alpha_i\chi_{A_i}, A_i\subset R_i$.
\medskip

First, set $A_1=R_1$ and $\varphi_1=\alpha_1\chi_{R_1}$ such that
$\int\varphi_1=\int_{6Q_i}f\omega_1$. Assume that $\varphi_1,
\ldots, \varphi_{k-1}$ have been constructed satisfying (\ref{cz4})
and
\begin{equation*}
\sum_{i=1}^{K-1}\varphi_i\leq \kappa \lambda,
\end{equation*}
where $\kappa$ is some constant which will be fixed later. There are two
cases:\\

\textbf{Case 1:} There exists some $i\in \{1,\ldots, k-1\}$ such
that $R_i\cap R_k\neq \emptyset$. Let $R_{s_1},\ldots, R_{s_m}$ be
the family of $R_1, \ldots, R_{k-1}$ such that $R_{s_j}\cap R_k\neq
\emptyset$. Since $l(R_{s_j})\leq l(R_k)$, $R_{s_j}\subset 3R_{k}$.
By using $R_k$ is $(3\times6^2, C_\lambda^{\log_2 3\times 6^2 +1})$-doubling
and (\ref{cz2}), we get
\begin{equation*}
\begin{aligned}
\sum_{j}|\varphi_{s_j}|&\leq \sum_{j}\int_{X}|f\omega_{s_j}|d\mu\\
&\leq C\sum_{j}\int_{X}\omega_{s_j}|f|d\mu\leq
C\sum_{j}\int_{3R_k}|f|d\mu\leq C\lambda \mu(3 .6^2R_k)\leq
C_1\lambda \mu(R_k).
\end{aligned}
\end{equation*}
Therefore,
\begin{equation*}
\mu\{\sum_j\|\varphi_{s_j}|>2C_1\lambda\}\leq \frac{\mu(R_k)}{2}.
\end{equation*}
Thus,
\begin{equation*}
\mu(A_k)\geq \frac{\mu(R_k)}{2}, A_k=R_k\cap
\{\sum_j|\varphi_{s_j}|\leq 2C_1\lambda\}.
\end{equation*}
The constant $\alpha_k$ will be chosen such that
$\int\varphi_k=\int_{Q_k}f\omega_kd\mu$ where
$\varphi_k=\alpha_k\chi_{A_k}$. Then we obtain
\begin{equation*}
\alpha_k\leq \frac{C}{\mu(A_k)}\int_{X}w_i|f|d\mu\leq
C\frac{2}{\mu(R_k)}\int_{\frac{1}{6^2}R_k}|f|d\mu\leq C_2\lambda \
\ (\text{by using (\ref{cz2})}).
\end{equation*}
If we choose $\kappa = 2C_1+C_2$, (\ref{cz5}) follows.\\

\textbf{Case 2:} $R_{i}\cap R_k =\emptyset$ for all $i=1,\ldots,
k-1$. Set $A_k=R_k$ and $\varphi_k=\alpha_k\chi_{R_k}$ such that
$\int\varphi_k=\int_{Q_k}f\omega_kd\mu$. We also get
(\ref{cz5}).\\

By the  construction of the functions $\varphi_i$, it is easy to
see that $\mu(R_i)\leq 2\mu(A_k)$. Hence,
\begin{equation*}
||\varphi_i||_{\infty}\mu(R_i)\leq C\alpha_i\mu(A_k)\leq
C\int_{X}|f\omega_i|d\mu.
\end{equation*}
When the collection of balls $\{Q_i\}$ is not finite, we can argue
as in \cite[p.134]{T1}. This completes the proofs of (a) and (b).\\

(c) Since $R_i$ is the smallest $(3\times6^2, C_\lambda^{\log_2 3\times 6^2
+1})$-doubling ball of the family $\{3\times 6^2 Q_i\}_{k\geq 1}$, one has
$K_{Q_i, R_i}\leq C$. For each $i$, we consider the atomic block
$b_i=fw_i-\varphi_i$ supported in ball $R_i$. By (\ref{cz1}) and
(\ref{cz6.1}) we have
$$
|b_i|_{\Hp}\leq \f{C}{\lambda^{p-1}}\int_X |fw_i|^pd\mu
$$
which implies
$$
|b|_{\Hp}\leq \f{C}{\lambda^{p-1}}\int_X \sum_i|fw_i|^pd\mu\leq
\f{C}{\lambda^{p-1}}\int_X
(\sum_iw_i)^p|f|^pd\mu=\f{C}{\lambda^{p-1}}\int_X |f|^pd\mu.
$$
Our proof is completed.\\

\medskip

Using the Calder\'on-Zygmund decomposition and a standard argument,
see for example \cite[pp.43-44]{J} (also \cite[p.135]{T1}), we
obtain the following interpolation result for a linear operator. For
clarity and completeness, we sketch the proof below.
\begin{thm}\label{interpolatioHardyRBMO}
Let $T$ be a linear operator which is bounded from $\HH$ into
$L^1(\mu)$ and from $L^\vc(\mu)$ into ${\rm RBMO}(\mu)$. Then $T$
can be extended to a bounded operator on $L^p(\mu)$ for all
$1<p<\vc$.
\end{thm}
\emph{Proof:} For simplicity we may assume that $\|\mu\|=\vc$. Let
$f$ be a function in $L^\vc(\mu)$ with bounded support satisfying
$\int fd\mu =0$. Let us recall that the set of all such functions is
dense in $L^p(\mu)$ for all $1<p<\vc$. For such functions $f$, we
need only to show that
\begin{equation}\label{eq1-InterpolationHardyBMO}
\|M^\sharp Tf\|_{L^p(\mu)}\leq C\|f\|_{L^p(\mu)}, \ \ 1<p<\vc.
\end{equation}
Once (\ref{eq1-InterpolationHardyBMO}) is proved,   Theorem \ref{interpolatioHardyRBMO}
follows from   Theorem \ref{sharpmaximalestimate2}.\\
For such a function $f$ and $\lambda>0$, we can decompose the function $f$
as in Theorem \ref{CZdecomposition}
$$
f:=b+g=\sum_{i}(w_if-\varphi_i)+g.
$$
By (\ref{cz3}) and (\ref{cz5}), we have $\|g\|_{L^\vc(\mu)}\leq C\lambda$, and by (\ref{cz8})
$$
\|b\|_{H^{1,p}_{at}(\mu)}\leq \f{C}{\lambda^{p-1}}\|f\|^p_{L^p(\mu)}.
$$
Since $T$ is bounded from $L^\vc(\mu)$ into RBMO$(\mu)$, we have
$$
\|M^\sharp Tg\|_{L^\vc(\mu)}\leq C_0\lambda.
$$
Therefore,
$$
\{M^\sharp Tf>(C_0+1)\lambda\}\subset \{M^\sharp Tb>\lambda\}.
$$
The fact that $M^\sharp$ is of weak type $(1,1)$ gives
$$
\mu\{M^\sharp Tb>\lambda\}\leq C\f{\|Tb\|_{L^1(\mu)}}{\lambda}.
$$
Moreover, since $T$ is bounded from $H^{1,\vc}_{at}(\mu)$ into $L^1(\mu)$, we have
$$
\|Tb\|_{L^1(\mu)}\leq C\|b\|_{H^{1,\vc}_{at}(\mu)}\leq \f{C}{\lambda^{p-1}}\|f\|^p_{L^p(\mu)}.
$$
This implies
$$
\mu\{M^\sharp Tf>(C_0+1)\lambda\}\leq C\f{\|f\|^p_{L^p(\mu)}}{\lambda^{p}}.
$$
So the sublinear operator $M^\sharp T$ is of weak type $(p,p)$ for all $1<p<\vc$.
By Marcinkiewicz interpolation theorem the operator $M^\sharp T$
is bounded for all $1<p<\vc$. This completes our proof.

\subsection{The weak $(1,1)$ boundedness of Calder\'on-Zygmund operators}

\begin{thm}\label{weak1-1}
If a Calder\'on-Zygmund operator $T$ is bounded on $L^2(\mu)$, then
$T$ is of weak type $(1,1)$.
\end{thm}
\emph{Proof:}
Let $f\in L^1(\mu)$ and $\lambda>0$. We can assume
that $\lambda>\beta_0 \|f\|_{L^1(\mu)}/\|\mu\|$. Otherwise, there is
nothing to prove. Using the same notations as in Theorem
\ref{CZdecomposition} with $R_i$ which is chosen as the smallest
$(3\times 6^2, C_\lambda^{\log_2 3\times 6^2 +1})$-doubling ball of
the family $\{3\times 6^2 Q_i\}_{k\geq 1}$,  we can write $f=g+b$,
with
$$
g=f\chi_{{}_{X\backslash \cup_{i}6Q_i}}+\sum_i \varphi_i
$$
and
$$
b:=\sum_{i}b_i=\sum_{i}(w_if-\varphi_i).
$$
Taking into account (\ref{cz1}), one has
$$
\mu(\cup_i 6^2Q_i)\leq\f{C}{\lambda}\int_{Q_i}|f|d\mu
\leq\f{C}{\lambda}\int_{X}|f|d\mu
$$
where in the last inequality we use the pairwise disjoint property
of family $\{Q_i\}_{i}$.\\
We need only to show that
$$
\mu\{x\in X\backslash \cup_{i}6^2Q_i: |T f(x)|>\lambda\}\leq
\f{C}{\lambda}\int_{X}|f|d\mu.
$$
We have
\begin{equation*}
\begin{aligned}
\mu\{x\in X\backslash \cup_{i}6^2Q_i: |T f(x)|>\lambda\}&\leq
\mu\{x\in X\backslash \cup_{i}6^2Q_i:
|T g(x)|>\lambda/2\}\\
&~~~+\mu\{x\in X\backslash \cup_{i}6^2Q_i: |T
b(x)|>\lambda/2\}:=I_1+I_2.
\end{aligned}
\end{equation*}
Let us estimate the term $I_1$ related to the ``good part'' first.
Since $|g|\leq C\lambda$ then
\begin{equation*}
\begin{aligned}
\mu\{x\in X\backslash \cup_{i}6^2Q_i: |T g(x)|>\lambda/2\}\leq
\f{C}{\lambda^2}\int|g|^2d\mu\leq \f{C}{\lambda}\int|g|d\mu.
\end{aligned}
\end{equation*}
Furthermore, we have
\begin{equation*}
\begin{aligned}
\int|g|d\mu&\leq \int_{X\backslash
\cup_{i}6Q_i}|f|d\mu+\sum_{i}\int_{R_i}|\varphi_i|\\
&\leq \int_{X}|f|d\mu+\sum_{i}\mu(R_i)\|\varphi_i\|_{L^\vc(\mu)}\\
&\leq \int_{X}|f|d\mu+C\sum_{i}\int_{X}|fw_i|d\mu\\
&\leq C\int_{X}|f|d\mu \ .
\end{aligned}
\end{equation*}
Therefore,
$$
\mu\{x\in X\backslash \cup_{i}6^2Q_i: |T g(x)|>\lambda/2\}\leq
\f{C}{\lambda}\int|f|d\mu.
$$
For the term $I_2$, we have
\begin{equation*}
\begin{aligned}
I_2&\leq \f{C}{\lambda}\sum_{i}\Big(\int_{X\backslash 2R_i}|T
b_i|d\mu+\int_{2R_i}|T \varphi_i|d\mu+\int_{2R_i\backslash 6^2Q_i}|T
w_if|d\mu\Big)\\
&\leq \f{C}{\lambda}\sum_{i}\Big(K_{i1}+K_{i2}+K_{i3}\Big)
\end{aligned}
\end{equation*}
 Note that $\int b_i d\mu =0$ for all $i$.
We have, by (\ref{cond2ofC-Z}),
\begin{equation*}
\begin{aligned}
K_{i1}=\int_{X\backslash 2R_i}|T b_i|d\mu &\leq C\int |b_i|d\mu\\
&\leq \int_{X}|fw_i|d\mu+\int_{R_i}|\varphi_i|
d\mu\\
&\leq \int_{X}|fw_i|d\mu+\mu(R_i)\|\varphi_i\|_{L^\vc(\mu)}\\
&\leq C\sum_{i}\int_{X}|fw_i|d\mu\\
&\leq C\sum_{i}\int_{X}|f|d\mu.
\end{aligned}
\end{equation*}
On the other hand, by the $L^2$ boundedness of $T$ and $R_i$ is
 a $(3\times 6^2, C_\lambda^{\log_2 3\times 6^2+1})$-doubling ball, we get
\begin{equation*}
\begin{aligned}
K_{i2}&\leq \Big(\int_{2R_i}|T
\varphi_i|^2\Big)^{1/2}(\mu(2R_i))^{1/2}\\
&\leq \Big(\int_{2R_i}|\varphi_i|^2\Big)^{1/2}(\mu(2R_i))^{1/2}\\
&\leq C\|\varphi_i\|_{L^\vc(\mu)}\mu(2R_i)\\
&\leq C\int |w_if|d\mu.
\end{aligned}
\end{equation*}
Moreover, taking into account the fact that supp$w_if\subset 6Q_i$, for $x\in
2R_i\backslash 6^2Q_i$ we have, by Lemma \ref{lem1-CZ},
\begin{equation*}
\begin{aligned}
K_{i3}&\leq C\int_{2R_i\backslash 6^2Q_i}\f{1}{\lambda(x_{Q_i},
d(x,x_{Q_i}))}\dx\times \int_{X}|w_if|d\mu.
\end{aligned}
\end{equation*}
Hence we obtain
$$
I_2\leq \f{C}{\lambda}\sum_{i}\int_{X}|w_if|d\mu \leq
\f{C}{\lambda}\sum_{i}\int_{X}|f|d\mu
$$
and the proof is completed.

\subsection{Cotlar inequality}
We note that from the weak type $(1,1)$ estimate
of $T$,  we can obtain a Cotlar inequality on $T$. More
precisely, we have  the following result.

\begin{thm}\label{Cotlar1}
Assume that  $T$ is a Calder\'on-Zygmund operator and that $T$ is bounded on
$L^2(X, \mu)$. Then there exists a constant $C>0$ such that for any
bounded function $f$ with compact support and $x\in X$ we have
$$
T_* f(x)\leq C\Big(M_{6,\eta}(Tf)(x)+M_{(5)}f(x)\Big)
$$
where
$$
M_{p,\rho}f(x)=\sup_{Q\ni x}\Big(\f{1}{\mu(\rho Q)}\int_Q|f|^p\Big).
$$
\end{thm}
\emph{Proof:} For any $\ep>0$ and $x\in X$, let $Q_x$ be the biggest
$(6, \beta)$-doubling ball centered $x$ of the form $6^{-k}\epsilon,
k\geq 1$ and $\beta>6^n$. Assume that $Q_x=B(x, 6^{-k_0}\epsilon)$.
Then, we can break $f=f_1+f_2$, where $f_1=f\chi_{\f{6}{5}Q_x}$.
Obviously $T_*f_1(x)=0$. This follows that
\begin{equation*}
T_*f(x)\leq |Tf_2(x)|+\Big|\int_{d(x,y)\leq
\ep}K(x,y)f_2(y)\dy\Big|=I_1+I_2.
\end{equation*}
Let us estimate $I_1$ first. For any $z\in Q_x$, we have
\begin{equation}\label{eq1-cotlar}
|Tf_2(x)|\leq |Tf_2(x)-Tf_2(z)|+|Tf(z)|+|Tf_1(z)|.
\end{equation}
On the other hand, it follows from  (\ref{cond2ofC-Z}) that
\begin{equation*}
\begin{aligned}
|Tf_2(x)-Tf_2(z)|&\leq \int_{X\backslash
\f{6}{5}Q_x}|K(x,y)-K(z,y)||f(y)|\dy\\
&\leq C \int_{X\backslash
\f{6}{5}Q_x}\f{d(x,z)^{\delta}}{d(x,y)^{\delta}\lambda(x,d(x,y))}|f(y)|\dy\\
&\leq C \int_{X\backslash
Q_x}\f{r_{Q_x}^{n+\delta}}{d(x,y)^{\delta}\lambda(x,d(x,y))}|f(y)|\dy\\
&\leq C \sum_{k=0}^\vc\int_{6^{k+1}Q_x\backslash 6^kQ_x}\f{r_{Q_x}^{\delta}}{d(x,y)^{\delta}\lambda(x,d(x,y))}|f(y)|\dy\\
&\leq C \sum_{k=0}^\vc\int_{6^{k+1}Q_x\backslash
6^kQ_x}\f{r_{Q_x}^{\delta}}{(6^kr_{Q_x})^{\delta}\lambda(x,6^kr_{Q_x})}|f(y)|\dy\\
&\leq C \sum_{k=0}^\vc 6^{-k}\f{\mu(6\times
6^{k+1}Q_x)}{\lambda(x,6^kr_{Q_x})}\f{1}{\mu(6\times
6^{k+1}Q_x)}\int_{6^{k+1}Q_x}|f(y)|\dy\\
&\leq C \sum_{k=0}^\vc 6^{-k}M_{(6)}f(x)=CM_{(6)}f(x).
\end{aligned}
\end{equation*}
This together with (\ref{eq1-cotlar}) implies
\begin{equation}\label{eq1-cotlar}
|Tf_2(x)|\leq CM_{(6)}f(x)+|Tf_1(z)|+ |Tf(z)|
\end{equation}
for all $z\in Q_x$.\\
At this stage, taking the $L^\eta(Q_x,\f{\dx}{\mu(Q)})$-norm with
respect to $z$, we have
$$
|Tf_2(x)|\leq
CM_{(6)}f(x)+\Big(\f{1}{\mu(Q_x)}\int_{Q_x}|Tf_1(z)|^\eta\dz\Big)^{1/\eta}+\Big(\f{1}{\mu(Q_x)}\int_{Q_x}|Tf(z)|^\eta\dz\Big)^{1/\eta}.
$$
By the Kolmogorov inequality and the weak type $(1,1)$ boundedness
of $T$, we have, for $\eta<1$,
\begin{equation*}
\begin{aligned}
\Big(\f{1}{\mu(Q_x)}\int_{Q_x}|Tf_1(z)|^\eta\dz\Big)^{1/\eta}&\leq
\f{1}{\mu(Q_x)}\int_{\f{6}{5}Q_x}|f_1(z)|\dz\\
&\leq \f{C}{\mu(15Q_x)}\int_{\f{6}{5}Q_x}|f_1(z)|\dz \ \ \ \text{(since $Q_x$ is $(6,\beta)$- doubling)}\\
&\leq CM_{(5)}f(x).
\end{aligned}
\end{equation*}
Furthermore, since $Q_x$ is $(6,\beta)$-doubling,
$$
\Big(\f{1}{\mu(Q_x)}\int_{Q_x}|Tf(z)|^\eta\dz\Big)^{1/\eta}\leq
CM_{\eta, 6}M(Tf)(x).
$$
Therefore, $I_1\leq CM_{(6)}f(x)+CM_{\eta, 6}M(Tf)(x)$.\\

For the term $I_2$ we have
\begin{equation*}
\begin{aligned}
I_2& \leq \int_{d(x,y)\leq \ep}|K(x,y)||f_2(y)|\dy\\
& \leq C\int_{B(x,\ep)\backslash B(x, 6^{-k_0}\ep)}\f{1}{\lambda(x,d(x,y))}|f_2(y)|\dy\\
& \leq C\sum_{k=0}^{k_0-1}\int_{B(x,6^{k+1-k_0}\ep)\backslash B(x, 6^{k-k_0}\ep)}\f{1}{\lambda(x,d(x,y))}|f_2(y)|\dy\\
& \leq C\sum_{k=0}^{k_0-1}\int_{B(x,6^{k+1-k_0}\ep)}\f{1}{\lambda(x,6^{k-k_0}\ep)}|f_2(y)|\dy\\
& \leq C\sum_{k=0}^{k_0-1}\f{\mu(x,6\times
6^{k+1-k_0}\ep)}{\lambda(x, 6\times
6^{k+1-k_0}\ep)}\f{1}{\mu(x,6\times
6^{k+1-k_0}\ep)}\int_{B(x,6^{k+1-k_0}\ep)}|f_2(y)|\dy\\
& \leq C\sum_{k=0}^{k_0-1}\f{\mu(x,6\times
6^{k+1-k_0}\ep)}{\lambda(x, 6\times 6^{k+1-k_0}\ep)}M_{(6)}f(x).
\end{aligned}
\end{equation*}
At this stage, by repeating the argument in the proof of Lemma
\ref{lem1-CZ}, we have
$$
\sum_{k=0}^{k_0-1}\f{\mu(x,6\times 6^{k+1-k_0}\ep)}{\lambda(x,
6\times 6^{k+1-k_0}\ep)}\leq C.
$$
Therefore,
$$
 I_2\leq CM_{(6)}f(x).
$$
This completes our proof.\\
\begin{rem}
From the boundedness of $M_{6,\eta}(\cdot)$ and $M_{(5)}(\cdot)$,
the Cotlar inequality tells us that if $T$ is bounded on $L^2(X,
\mu)$ then the maximal operator $T_*$ is bounded on $L^p(X, \mu)$
for $1<p<\vc$. Note that the endpoint estimate of $T_*$ will be
investigated in \cite{AD}.

The Calder\'on-Zygmund decomposition Theorem \ref{CZdecomposition}
does not require the property (v) of $\lambda(\cdot, \cdot)$.
\end{rem}
\section{The boundedness of Calder\'on-Zygmund operators}

The main results of this section are Theorems \ref{boundedofCZOonRBMO},
\ref{boundedofCZOonHardyspace} and \ref{boundednessofcommutators}.

\subsection{The boundedness of Calder\'on-Zygmund operators from $L^{\infty}$ to  RBMO
space}

The following result shows that on a non-homogeneous space $(X, \mu)$, a Calder\'on-Zygmund operator
which is bounded on $L^2$ is also bounded from $L^\vc(\mu)$ into the regularized BMO
space RBMO$(\mu)$.

\begin{thm}\label{boundedofCZOonRBMO}
Assume that $T$ is a Calder\'on-Zygmund operator and $T$ is bounded on $L^2(\mu)$, then $T$
is bounded from $L^\vc(\mu)$ into RBMO$(\mu)$. Therefore, by
interpolation and duality, $T$ is bounded on $L^p(\mu)$ for all
$ 1 <  p<\vc$.
\end{thm}
\emph{Proof:} We use the RBMO characterizations (\ref{cond1-RBMOdefn3})
and (\ref{cond2-RBMOdefn3}). The condition (\ref{cond1-RBMOdefn3})
can be obtained by the standard method used in the case of doubling
measure. We omit the details here.\\

We will check (\ref{cond2-RBMOdefn3}). To do this, we have to show
that
$$
|m_Q(T f)-m_R(T f)|\leq
CK_{Q,R}\Big(\f{\mu(6Q)}{\mu(Q)}+\f{\mu(6R)}{\mu(R)}\Big)\|f\|_{L^\vc(\mu)}
$$
for all $Q\subset R$.\\
Let $N$ be the first integer $k$ such that $R\subset 6^kQ$. We
denote $Q_R=6^{N+1}Q$. Then for $x\in Q$ and $y\in R$, we set
\begin{equation*}
\begin{aligned}
T f(x)-T f(y)&=T f\chi_{6Q}(x)+T f\chi_{6^{N}Q\backslash 6Q}(x)+T
f\chi_{X\backslash
Q_R}(x)\\
&~~~ -(T f\chi_{Q_R}(y)+ T f\chi_{X\backslash Q_R}(y))\\
&\leq |T f\chi_{6Q}(x)|+|T
f\chi_{6^{N}Q\backslash 6Q}(x)|\\
&~~~+|T f\chi_{X\backslash Q_R}(x)- T f\chi_{X\backslash
Q_R}(y)|+|T f\chi_{Q_R}(y)|\\
&\leq I_1+I_2+I_3+I_4.
\end{aligned}
\end{equation*}
Let us estimate $I_3$ first. We have
\begin{equation*}
\begin{aligned}
I_3&\leq \int\limits_{X\backslash Q_R}|K(x,z)-K(y,z)||f(z)|\dz\\
&\leq \sum_{k=N+1}^\vc \int\limits_{6^{k+1}Q\backslash 6^kQ}\f{d(x,y)^\de}{d(x,z)^\de\lambda(x,d(x,z))}|f(z)|\dz\\
&\leq \sum_{k=N+1}^\vc
6^{-(k-N)\delta}\f{\mu(6^{k+1}Q)}{\la(x,6^{k-1}r_Q)}\|f\|_{L^\vc(\mu)}\\
&\leq C\sum_{k=N+1}^\vc
6^{-(k-N)\delta}\f{\mu(6^{k+1}Q)}{\la(x,6^{k+1}r_Q)}\|f\|_{L^\vc(\mu)}\\
&\leq C\sum_{k=N+1}^\vc
6^{-(k-N)\delta}\|f\|_{L^\vc(\mu)}=C\|f\|_{L^\vc(\mu)},
\end{aligned}
\end{equation*}
where in the last inequality we use the fact that $\mu(6^{k+1}Q)\leq {\la(x,6^{k+1}r_Q)}$, since $x\in Q\subset 2^{k+1}Q$.\\
As to the term $I_2$, we have
\begin{equation}\label{eq1-RBMOboundedproof}
\begin{aligned}
T f\chi_{6^{N}Q\backslash 6Q}(x)&\leq
\int_{6^{N}Q\backslash 6Q}|K(x,y)||f(y)|\dy\\
&\leq \int_{6^{N}Q\backslash 6Q}\f{C}{\lambda(x_Q,d(x_Q,y))}|f(y)|\dy\\
&\leq K_{6Q,6^{N}Q}\|f\|_{L^\vc(\mu)}.
\end{aligned}
\end{equation}
Therefore, $I_2\leq CK_{Q,R}\|f\|_{\vc}$. So, we get
\begin{equation*}
\begin{aligned}
T f(x)-T f(y)&=T f\chi_{6Q}(x)+CK_{Q,R}\|f\|_{L^\vc(\mu)}+|T
f\chi_{Q_R}(y)|+ C\|f\|_{L^\vc(\mu)}.
\end{aligned}
\end{equation*}
Taking the mean over $Q$ and $R$ for $x$ and $y$, respectively, we
have
\begin{equation*}
\begin{aligned}
|m_Q(T f)-m_R(T f)|&\leq m_Q(|T
f\chi_{6Q}|)+CK_{Q,R}\|f\|_{L^\vc(\mu)}+|T f\chi_{Q_R}(y)|\\
&+ C\|f\|_{L^\vc(\mu)} +m_R(T f\chi_{Q_R}).
\end{aligned}
\end{equation*}
For the boundedness on $L^2(\mu)$ of $T$, we have
\begin{equation*}
\begin{aligned}
m_Q(|T f\chi_{6Q}|)&\leq \Big(\f{1}{\mu(Q)}\int_Q|T
f\chi_{6Q}|^2\Big)^{1/2}\\
&\leq C\Big(\f{\mu(6Q)}{\mu(Q)}\Big)^{1/2}\|f\|_{L^\vc(\mu)}\\
&\leq C\Big(\f{\mu(6Q)}{\mu(Q)}\Big)\|f\|_{L^\vc(\mu)}.
\end{aligned}
\end{equation*}
Next, we write
$$
m_R(T f\chi_{Q_R})\leq m_R(T f\chi_{Q_R\cap 6R})+m_R(T f\chi_{Q_R
\backslash 6R}).
$$
By similar argument in estimate of $m_Q(|T f\chi_{6Q}|)$, the term
$m_R(T f\chi_{Q_R\cap 6R})$ is dominated by
$$
C\Big(\f{\mu(6R)}{\mu(R)}\Big)\|f\|_{L^\vc(\mu)}.
$$
The second term $m_R(T f\chi_{Q_R \backslash 6R})$ can be treated as
in (\ref{eq1-RBMOboundedproof}). Since $r_{Q_R}\approx r_R$, we have
$$
m_R(T f\chi_{Q_R \backslash 6R})\leq C\|f\|_{L^\vc(\mu)}
$$
To sum up, we have
\begin{equation*}
\begin{aligned}
|m_Q(T f)-m_R(T f)|&\leq
CK_{Q,R}\|f\|_{L^\vc(\mu)}+\Big(\f{\mu(6Q)}{\mu(Q)}+\f{\mu(6R)}{\mu(R)}\Big)\|f\|_{L^\vc(\mu)}\\
&\leq
CK_{Q,R}\Big(\f{\mu(6Q)}{\mu(Q)}+\f{\mu(6R)}{\mu(R)}\Big)\|f\|_{L^\vc(\mu)}.
\end{aligned}
\end{equation*}

\begin{rem}
By similar argument in \cite[Theorem 2.11]{T1}, we can replace the
assumption of $L^2(\mu)$ boundedness by the weaker assumption: for
any ball $B$ and any function $a$ supported on $B$,
$$
\int_B |T a|d\mu\leq C\|a\|_{L^\vc}\mu(6B)
$$
uniformly on $\epsilon > 0$.
\end{rem}

\subsection{The boundedness of Calder\'on-Zygmund operators on Hardy spaces}

We now show that an $L^2$ bounded Calder\'on-Zygmund operator maps the atomic
Hardy space boundedly into $L^1$.

\begin{thm}\label{boundedofCZOonHardyspace}
Assume that $T$ is a  Calder\'on-Zygmund operator and $T$ is bounded on $L^2(X, \mu)$, then
$T$ is bounded from $\HH $ into $L^1(X, \mu)$. Therefore, by
interpolation and duality, $T$ is bounded on $L^p(\mu)$ for all
$1<p < \infty$.

\end{thm}
\emph{Proof:} By \cite[Lemm 4.1]{HoM}, it is enough to show that
\begin{equation}\label{eq1-CZHardy}
\|T b\|_{L^1(\mu)} \leq C|b|_{\HH}
\end{equation}
for any atomic block $b$ with supp$b\subset B$ and
$=\sum_{j}\lambda_j a_j$ where the $a_j$'s are functions satisfying
(a) and (b) in definition of atomic blocks. At this stage we can
use the same argument as in \cite[Theorem 4.2]{T1} with minor
modifications as in Theorem \ref{boundedofCZOonRBMO} to obtain the
estimate (\ref{eq1-CZHardy}). We omit the
details here.\\

\subsection{Commutators of Calder\'on-Zygmund operators with RBMO functions}
In this section we assume that the dominating function $\lambda$
satisfies $\lambda(x,ar)=a^m\lambda(x,r)$ for all $x\in X$ and $a,
r>0$. Then, for two balls $B$, $Q$ such that $B\subset Q$ we can define the coefficient
$K'_{B, Q}$ as follows: let $N_{B, Q}$ be the smallest integer
satisfying $6^{N_{B, Q}}r_B\geq r_{Q}$, we set
\begin{equation}\label{coefficientK1}
K'_{B,Q}:=1+\sum_{k=1}^{N_{B,Q}}\f{\mu(6^kB)}{\lambda(x_B, 6^kr_B)}.
\end{equation}
It is not difficult to show that the coefficient $K_{B, Q}\approx
K'_{B, Q}$. Note that in the definition of $K'_{B, Q}$ we can replace
$6$ by any number $\eta>1$.\\

 To establish the boundedness of
commutators of Calder\'on-Zygmund operators with RBMO functions, we
need the following two lemmas. Note that these lemmas are similar to
those in \cite{T1}. However, due to the difference of choices of
coefficient $K_{Q,R}$, we would like to provide the proof for the
first one. Meanwhile, the proof of Lemma \ref{lem2-commutator} is
completely analogous to that of Lemma 9.3 in \cite{T1}, hence we
omit the details.
\begin{lem}\label{lem1-commutator}
If $B_i=B(x_0,r_i), i=1,\ldots, m$ are concentric balls $B_1\subset
B_2\subset\ldots\subset B_m$ with $K_{B_i,B_{i+1}}>2$ for
$i=1,\ldots, m-1$ then
\begin{equation}
\sum_{i=1}^{m-1}K_{B_i,B_{i+1}}\leq 2K_{B_1,B_{m}}.
\end{equation}
\end{lem}
\emph{Proof:} By definition,
$$
K_{B_i,B_{i+1}}=1+\int\limits_{r_i\leq d(x,x_0)\leq
r_{i+1}}\f{1}{\la(x_0,d(x,x_0))}\dx.
$$
Since $K_{B_i,B_{i+1}}>2$, we have
$$
K_{B_i,B_{i+1}}<2\int\limits_{r_i\leq d(x,x_0)\leq
r_{i+1}}\f{1}{\la(x_0,d(x,x_0))}\dx
$$
for all $i=1,\ldots, m-1$.\\
This implies
\begin{equation*}
\begin{aligned}
\sum_{i=1}^{m-1}K_{B_i,B_{i+1}}&<2\sum_{i=1}^{m-1}\int\limits_{r_i\leq
d(x,x_0)\leq r_{i+1}}\f{1}{\la(x_0,d(x,x_0))}\dx\\
&\leq 2 \int\limits_{r_1\leq
d(x,x_0)\leq r_{m}}\f{1}{\la(x_0,d(x,x_0))}\dx\\
&\leq 2K_{B_1,B_{m}}.
\end{aligned}
\end{equation*}
\begin{lem}\label{lem2-commutator}
There exists some constant $P_0$ such that if $x\in X$ is some fixed
point and $\{f_B\}_{B\ni x}$ is collection of numbers such that
$|f_Q-f_R|\leq C_x$ for all doubling balls $Q\subset R$ with $x\in
Q$ and $K_{Q,R}\leq P_0$, then
$$
|f_Q-f_R|\leq CK_{Q,R}C_x \ \ \text{for all doubling balls $Q\subset
R$ with $x\in Q$}.
$$
\end{lem}
\begin{thm}\label{boundednessofcommutators}
If $b\in {\rm RBMO}(\mu)$ and $T$ is a Calder\'on-Zygmund bounded on
$L^2(\mu)$, then the commutator $[b, T]$ defined by
$$
[b,T](f)=bT(f)-T(bf)
$$
is bounded on $L^p(\mu)$ for $1<p<\vc$.
\end{thm}
\emph{Proof:} For $1<p<\vc$ we will show that
\begin{equation}\label{pointwiseestforcommutators}
\MS([b, T]f)(x)\leq C\|b\|_{{\rm
RBMO}(\mu)}\Big(M_{p,5}f(x)+M_{p,6}Tf(x)+T_*f(x)\Big),
\end{equation}
where
$$
M_{p,\rho}f(x)=\sup_{Q\ni x}\Big(\f{1}{\mu(\rho Q)}\int_Q|f|^p\Big).
$$
Once (\ref{pointwiseestforcommutators}) is proved, it follows from the
boundedness of $T_*$ on $L^p(\mu)$ and $M_{p,\rho}$ on
$L^r(\mu),r>p$ and $\rho\geq 5$, and from a standard argument that we can obtain
the boundedness
of $[b, T]$ on $L^p(\mu)$.\\

\noindent Let $\{b_B\}$ be a family of numbers satisfying
$$
\int_B|b-b_B|d\mu\leq 2\mu(6B)\|b\|_{{\rm RBMO}}
$$
for balls $B$, and
$$
|b_Q-b_R|\leq 2K_{Q,R}\|b\|_{{\rm RBMO}}
$$
for balls $Q\subset R$.
Denote
$$
h_Q:=m_Q(T((b-b_Q)f\chi_{X\backslash\f{6}{5}Q}).
$$
We will show that
\begin{equation}\label{eq1-commutator}
\f{1}{\mu(6B)}\int_B|[b,T]f-h_Q|d\mu\leq C\|b\|_{{\rm
RBMO}}(M_{p,5}f(x)+M_{p,6}Tf(x))
\end{equation}
for all $x$ and $B$ with $x\in B$, and
\begin{equation}\label{eq2-commutator}
|h_Q-h_R|\leq C\|b\|_{{\rm RBMO}}(M_{p,5}f(x)+T_*f(x))K^2_{Q,R}
\end{equation}
for all  $x\in Q\subset R$.\\
The proof of (\ref{eq1-commutator}) is similar to that in
Theorem 9.1 in \cite{T1} with minor modifications and we omit
it here. \\
It remains to check (\ref{eq2-commutator}). For two balls $Q\subset
R$, let $N$ be an integer  such that $(N-1)$ is the smallest
number satisfying $r_R\leq 6^{N-1}r_Q$. Then, we break the term
$|h_Q-h_R|$ into five terms:
\begin{equation*}
\begin{aligned}
|m_Q(T((b-b_Q)&f\chi_{X\backslash\f{6}{5}Q})-m_R(T((b-b_R)f\chi_{X\backslash\f{6}{5}R})|\\
&\leq
|m_Q(T((b-b_Q)f\chi_{6Q\backslash\f{6}{5}Q})|+|m_Q(T((b_Q-b_R)f\chi_{X\backslash
6Q})|\\
&~~~+|m_Q(T((b-b_R)f\chi_{6^NQ\backslash 6Q})|\\
&~~~+|m_Q(T((b-b_R)f\chi_{X\backslash
6^NQ})-m_R(T((b-b_R)f\chi_{X\backslash 6^NQ})|\\
&~~~+|m_R(T((b-b_R)f\chi_{6^NQ\backslash \f{6}{5}R})\\
&=M_1+M_2+M_3+M_4+M_5.
\end{aligned}
\end{equation*}
Let us estimate $M_1$ first. For $y\in Q$ we have, by Proposition
3.2
\begin{equation*}
\begin{aligned}
|T((b-b_Q)&f\chi_{6Q\backslash\f{6}{5}Q})(x)|\\
&\leq
\f{C}{\la(x,r_Q)}\int_{6Q}|b-b_Q||f|d\mu\\
&\leq
\f{\mu(30Q)}{\la(x, 30r_Q)}\Big(\f{1}{\mu(5\times 6Q)}\int_{6Q}|b-b_Q|^{p'}d\mu\Big)^{1/p'}\Big(\f{1}{\mu(5\times 6Q)}\int_{6Q}|f|^{p}d\mu\Big)^{1/p}\\
&C\|b\|_{{\rm RBMO}}M_{p,5}f(x).
\end{aligned}
\end{equation*}
The term $M_5$ can be treated by similar way. So, we have
$M_1+M_5\leq C\|b\|_{{\rm RBMO}}M_{p,5}f(x).$ For the term $M_2$, we
have for $x,y\in Q$
\begin{equation*}
\begin{aligned}
|Tf\chi_{X\backslash 6Q}(y)|&= \Big|\int_{X\backslash 6Q}K(y,z)f(z)d\mu(z)\Big|\\
&\leq \Big|\int_{X\backslash 6Q}K(x,z)f(z)d\mu(z)\Big| +\int_{X\backslash 6Q}|K(y,z)-K(x,z)||f(z)|d\mu(z)\\
& \leq T_*f(x)+CM_{p,5}f(x).
\end{aligned}
\end{equation*}
This implies
$$
|m_Q(T((b_Q-b_R)f\chi_{X\backslash 6Q})|\leq
CK_{Q,R}(T_*f(x)+M_{p,6}f(x)).
$$
For the term $M_4$, we have, for $y,z\in R$
\begin{equation*}
\begin{aligned}
|T((b-b_R)&f\chi_{X\backslash 6^NQ}(y)-T((b-b_R)f\chi_{X\backslash
6^NQ}(z)|\\
&\leq \int_{X\backslash 2R}|K(y,w)-K(z,w)||(b(w)-b_R)||f(x)|\dw\\
&\leq \int_{X\backslash 2R}\f{d(y,z)^{\de}}{d(w,y)^\delta\la(y,d(w,y))}|(b(w)-b_R)||f(x)|\dw\\
&\leq \sum_{k=1}^\vc\int_{2^{k+1}R\backslash 2^{k}R}\f{d(y,z)^{\de}}{d(w,y)^\delta\la(y,d(w,y))}|(b(w)-b_R)||f(x)|\dw\\
&\leq \sum_{k=1}^\vc 2^{-k\de}\f{1}{\la(y, 2^{k-1}r_R)}\int_{2^{k+1}R}|(b(w)-b_R)||f(x)|\dx\\
\end{aligned}
\end{equation*}
By H\"older inequality we have
\begin{equation*}
\begin{aligned}
|T((b-b_R)&f\chi_{X\backslash 6^NQ}(y)-T((b-b_R)f\chi_{X\backslash
6^NQ}(z)|\\
&\leq C\sum_{k=1}^\vc 2^{-k\de}\f{\mu(5\times 2^{k+1}R)}{\la(y, 5\times 2^{k+1}r_R)}\Big(\f{1}{\mu(5\times 2^{k+1}R)}\int_{2^{k+1}R}|b-b_R|^{p'}d\mu\Big)^{1/p'}\\
&~~~~\Big(\f{1}{\mu(5\times 2^{k+1}R)}\int_{2^{k+1}R}|f|^{p}d\mu\Big)^{1/p}\\
&\leq C\sum_{k=1}^\vc 2^{-k\de}\Big[\Big(\f{1}{\mu(5\times 2^{k+1}R)}\int_{2^{k+1}R}|b-b_{2^{k+1}R}|^{p'}d\mu\Big)^{1/p'}\\
&~~~~+\Big(\f{1}{\mu(5\times 2^{k+1}R)}\int_{2^{k+1}R}|b_R-b_{2^{k+1}R}|^{p'}d\mu\Big)^{1/p'}\Big]\Big(\f{1}{\mu(5\times 2^{k+1}R)}\int_{2^{k+1}R}|f|^{p}d\mu\Big)^{1/p}\\
&\leq C\sum_{k=1}^{\vc}C(k+1)2^{-k\de}\|b\|_{{\rm
RBMO}}M_{p,5}f(x)\\
&\leq C\|b\|_{{\rm RBMO}}M_{p,5}f(x).
\end{aligned}
\end{equation*}
Taking the mean over $Q$ and $R$ for $y$ and $z$ respectively, we
obtain
$$
M_4\leq C\|b\|_{{\rm RBMO}}M_{p,5}f(x).
$$
Concerning the last estimate for $M_3$, we have for $y\in Q$
\begin{equation}\label{eq3-commutator}
\begin{aligned}
|T((b-b_R)&f\chi_{6^NQ\backslash 6Q}(y)|\\& \leq
C\sum_{k=1}^{N-1}\f{1}{\la(6^kQ)}\int_{6^{k+1Q}\backslash
6^kQ}|b-b_R||f|d\mu\\
&\leq C\sum_{k=1}^{N-1}\f{1}{\la(y,
6^kQ)}\Big[\int\limits_{6^{k+1}Q\backslash
6^kQ}|b-b_{5^{k+1}Q}||f|d\mu+\int\limits_{6^{k+1Q}\backslash
6^kQ}|b_R-b_{6^{k+1}Q}||f|d\mu\Big]\\
&\leq C\sum_{k=1}^{N-1}\f{\mu(5\times 6^{k+1}Q)}{\la(x_Q,
6^kQ)}\Big[\f{1}{\mu(6^{k+2}Q)}\int\limits_{6^{k+1}Q\backslash
6^kQ}|b-b_{6^{k+1}Q}||f|d\mu\\
&~~~~+\f{1}{\mu(5\times  6^{k+1}Q)}\int\limits_{6^{k+1Q}\backslash
6^kQ}|b_R-b_{6^{k+1}Q}||f|d\mu\Big] .\\
\end{aligned}
\end{equation}
By H\"older inequality and a similar argument to the estimate of the term
$M_4$, we have
$$
\f{1}{\mu(5\times 6^{k+2}Q)}\int\limits_{6^{k+1}Q\backslash
6^kQ}|b-b_{6^{k+1}Q}||f|d\mu\leq \|b\|_{{\rm RBMO}}M_{p,5}f(x)
$$
and
$$
\f{1}{\mu(5\times 6^{k+1}Q)}\int\limits_{6^{k+1Q}\backslash
6^kQ}|b_R-b_{6^{k+1}Q}||f|d\mu\leq CK_{Q,R}\|b\|_{{\rm
RBMO}}M_{p,5}f(x).
$$
These two estimates together with (\ref{eq3-commutator}) give
$$
|T((b-b_R)f\chi_{6^NQ\backslash 6Q}(y)|\leq CK^2_{Q,R}\|b\|_{{\rm
RBMO}}M_{p,5}f(x).
$$
This implies $M_3\leq CK^2_{Q,R}\|b\|_{{\rm RBMO}}M_{p,5}f(x).$ From
the estimates $M_1, M_2, M_3, M_4, M_5$, we obtain (\ref{eq2-commutator}).\\
To obtain (\ref{pointwiseestforcommutators}) from
(\ref{eq1-commutator}) and (\ref{eq2-commutator}), we use a trick of
\cite{T1}. From (\ref{eq1-commutator}), if $Q$ is a doubling ball
and $x\in Q$, we have
\begin{equation}\label{eq4-commutator}
\begin{aligned}
|m_Q([b,T]f)-h_Q|&\leq \f{1}{\mu(Q)}\int_Q |[b,T]f-h_Q|d\mu\\
&\leq C\|b\|_{{\rm RBMO}}(M_{p,5}f(x)+M_{p,6}Tf(x)).
\end{aligned}
\end{equation}
Also, for any ball $Q\ni x$ (non doubling, in general),
$K_{Q,\wt{Q}}\leq C$, and then by (\ref{eq1-commutator}) and
(\ref{eq2-commutator}) we have
\begin{equation}\label{eq5-commutator}
\begin{aligned}
\f{1}{\mu(6Q)}&\int_Q|[b,T]f-m_{\wt{Q}}[b,T]f|d\mu\\
&\leq \f{1}{\mu(6Q)}\int_Q|[b,T]f-h_Q|d\mu+|h_Q-h_{\wt{Q}}|+|h_{\wt{Q}}-m_{\wt{Q}}[b,T]f|\\
& \leq C\|b\|_{{\rm
RBMO}(\mu)}\Big(M_{p,5}f(x)+M_{p,6}Tf(x)+T_*f(x)\Big).
\end{aligned}
\end{equation}
In addition, for all doubling balls $Q\subset R$ with $x\in Q$ such
that $K_{Q,R}\leq P_0$ where $P_0$ is a constant in Lemma
\ref{lem2-commutator}, by (\ref{eq2-commutator}) we have
$$
|h_Q-h_R|\leq C \|b\|_{{\rm
RBMO}(\mu)}\Big(M_{p,5}f(x)+T_*f(x)\Big)P_0^2.
$$
Due to Lemma \ref{lem2-commutator} we get
$$
|h_Q-h_R|\leq C \|b\|_{{\rm
RBMO}(\mu)}\Big(M_{p,5}f(x)+T_*f(x)\Big)K_{Q,R},
$$
for all doubling balls $Q\subset R$ with $x\in Q$. At this stage,
applying (\ref{eq3-commutator}), we obtain
\begin{equation*}
\begin{aligned}
m_Q([b,T]f)&-m_{R}([b,T]f)\\
& \leq C\|b\|_{{\rm
RBMO}(\mu)}\Big(M_{p,5}f(x)+M_{p,6}Tf(x)+T_*f(x)\Big)K_{Q,R}.
\end{aligned}
\end{equation*}
This completes our proof.\\

\begin{rem}
  As mentioned  earlier in this paper, the  results
of this article  still hold when $X$ is a quasi-metric space. Indeed, one can see
 that the main problem in quasi-metric space setting is that
the covering lemma, Lemma 2.1, may not be true. However, instead of
using this covering property, we can adapt the covering lemma in
\cite[Lemma 3.1]{FGL} to our situation. This problem is not
difficult and we   leave it to the interested reader.
\end{rem}

\medskip
\noindent Department of Mathematics, Macquarie University, NSW 2109, Australia and \\
Department of Mathematics, University of Pedagogy, Ho chi Minh city, Vietnam \\
Email: the.bui@mq.ed.au

\medskip
\noindent Department of Mathematics, Macquarie University, NSW 2109, Australia \\
Email: xuan.duong@mq.edu.au


\begin{thebibliography}{99}

\bibitem[AD]{AD} T. A. Bui  and X. T.  Duong, Endpoint estimates for maximal operator and boundedness
of maximal commutators, preprint.

\bibitem[CW1]{CW1} R. Coifman and G. Weiss, {\it Analyse harmonique
non-commutative sur certains espaces homog\`enes}, Lecture Notes in
Mathematics {\bf 242}. Springer, Berlin-New York, 1971.

\bibitem[CW2]{CW2} R. Coifman and G. Weiss, Extensions of Hardy spaces and their use in
analysis, Bull. Amer. Math. Soc. Volume 83 (1977), 569-645.
\bibitem[D]{D} G. David, Completely uncrectifiable 1-sets on the
plane have vanishing analytic capacity, Pr\'epublications
Math\'ematiques d'Orsay, 1997 (61), 1-94.

\bibitem[FGL]{FGL} G. Di Fazio,
C. E. Guti\'errez, and E. Lanconelli, Covering theorems,
inequalities on metric spaces and applications to PDE's, Math. Ann.
341 (2008), no. 2, 255-291.

\bibitem[FS]{FS} G. Folland and E.M. Stein, {\it Hardy spaces on Homogeneous
Groups}, Princeton  Univ. Press, 1982.
\bibitem[J]{J}J.-L. Journ\'e, Calder\'on-Zygmund operators, pseudo-differential operators and the
Cauchy integral of Calder�on. Lecture Notes in Math. 994, Springer,
1983

\bibitem[He]{He} J. Heinonen, Lectures on analysis on metric spaces, Universitext, Springer-Verlag, New York, 2001.

\bibitem[HoM]{HoM} S. Hofmann and S. Mayboroda, Hardy and BMO spaces associated to divergence
form elliptic operators, Math. Ann. 344 (2009), no.1, 37-116.

\bibitem[Hy]{H} T. Hyt\"onen, A framework for non-homogenous analysis on metric spaces, and RBMO spaces of Tolsa, preprint.

\bibitem[HyM]{HM} T. Hyt\"onen and H. Martikainen, Non-homogeneous Tb theorem and random dyadic cubes on metric measure
spaces, preprint.


\bibitem[HyYY]{HYY} T. Hyt\"onen, Dachun Yang and Dongyong Yang, The Hardy Space $H^1$ on Non-homogeneous Metric
Spaces, preprint.

\bibitem[NTV1] {NTV1} F. Nazarov, S. Treil and A. Volberg,  Cauchy integral and Calder\'on-Zygmund
operators on nonhomogeneous spaces, Internat. Math. Res. Notices, Vol 15,
1997, p. 703 - 726.


\bibitem[NTV2] {NTV2} F. Nazarov, S. Treil and A. Volberg, Weak type estimates and Cotlar inequalities
for  Calder\'on-Zygmund
operators on nonhomogeneous spaces, Internat. Math. Res. Notices, Vol 9,
1998, p. 463 - 487.


\bibitem[NTV3]{NTV3}  F. Nazarov, S. Treil and A. Volberg, The $Tb$- theorem on non-homogeneous
spaces, Acta Math., Vol 190, 2003, No 2, p. 151 - 239.



\bibitem[MMNO]{MMNO} J. Mateu, P. Mattila, A. Nicolau, J. Orobitg, BMO for non
doubling measures, Duke Math. J., 102 (2000), 533-565.


\bibitem[St]{St} E.M. Stein,  {\it Harmonic analysis: Real variable
methods, orthogonality and oscillatory integrals}, Princeton Univ.
Press, Princeton, NJ, (1993).

\bibitem[Tch]{Tch} E. Tchoundja, Carleson measures for the generalized Bergman
spaces via a T (1)-type theorem, Ark. Mat. 46 (2008), no. 2,
377-406.
\bibitem[T1]{T1} X. Tolsa, BMO, $H^1$, and Calder\'on-Zygmund operators
for non doubling measures, Math. Ann. 319 (2001), 89-149.

\bibitem[T2]{T2} X. Tolsa, A proof of the weak (1,1) inequality for singular integrals with non
doubling measures based on a Calder\'on-Zygmund decomposition, Publ.
Mat. 45 (2001), 163-174

\bibitem[V]{V} J. Verdera. On the $T(1)$ theorem for the Cauchy integral. Ark.
Mat. 38 (2000), 183-199

\bibitem[VW]{VW} A. Volberg and B.
D. Wick, Bergman-type singular operators and the characterization of
Carleson measures for Besov-Sobolev spaces on the complex ball,
preprint (2009), arXiv:0910.1142.
\end{thebibliography}
\end{document}